\documentstyle[12pt]{amsart}
\newcommand{\nn}[1]{(\ref{#1})}

\newtheorem{theorem}{Theorem}[section]
\newtheorem{lemma}[theorem]{Lemma}
\newtheorem{proposition}[theorem]{Proposition}
\newtheorem{definition}[theorem]{Definition}

\newtheorem{corollary}[theorem]{Corollary}

\newtheorem{remark}[theorem]{Remark}

\newcommand{\sfrac}[2]{{\textstyle \frac{#1}{#2}}}

\newcommand{\toiso}{\usebox{\ttoiso}}
\newsavebox{\ttoiso}
\sbox{\ttoiso}{\begin{picture}(25,12)(-5,0)
\put(0,3){$\simeq$}
\put(0,-3){$\longrightarrow$}
\end{picture}}

\newcommand{\br}{\mbox{$\Bbb R$}}
\newcommand{\C}{\mbox{$\Bbb C$}}
\newcommand{\Do}{\mbox{\sf D}}
\newcommand{\Dok}{{\mbox{\sf D}^{(k)}}}

\newcommand{\cd}{\partial}
\newcommand{\nd}{\nabla}

\newcommand{\Rho}{\mbox{\sf P}}
\newcommand{\up}{\Upsilon}

\newcommand{\ce}{{\cal E}}
\newcommand{\cf}{{\cal F}}
\newcommand{\cg}{{\cal G}}
\newcommand{\ch}{{\cal H}}

\newcommand{\vep}{{\varepsilon}}
\newcommand{\ep}{{\epsilon}}

\newcommand{\vol}{\mbox{\large $\epsilon$}}
\newcommand{\cs}{{\cal S}}

\newcommand{\cj}{{\cal J}}
\newcommand{\cjk}{{{\cal J}^{k}}}
\newcommand{\cu}{{\cal U}}
\newcommand{\cv}{{\cal V}}
\newcommand{\cw}{{\cal W}}

\newcommand{\vb}{{\vphantom{b}}}

\newcommand{\bbox}                        
{\mbox{$
\begin{picture}(8,8)(0,0)
\put(0,0){$\Box$}
\end{picture}$}}

\newcommand{\tbox}{\widetilde{\bbox}}

\newcommand{\mb}                         
{\mbox{$
\begin{picture}(4,4)(0,0)
\put(0,5){\line(1,0){5}}
\put(0,0){\line(1,0){5}}
\put(0,0){\line(0,1){5}}
\put(5,0){\line(0,1){5}}
\end{picture}$}}

\newcommand{\yiiI}{\mbox{%
$\begin{picture}(13,7)(-1,1)
\put(0,0){\mb}
\put(5,0){\mb}
\end{picture}$}}

\newcommand{\yiIIi}{\mbox{%
$\begin{picture}(12,12)(-1,1)
\put(0,5){\mb}
\put(5,5){\mb}
\put(0,0){\mb}
\end{picture}$}}

\newcommand{\weylsymm}
{
\mbox{$\begin{picture}(20,21)(-1,0)
\put(0,3){\tiny ${q-1}$}
\put(9,10){\vector(0,1){6}}
 \put(9,1){\vector(0,-1){6}}
\end{picture}$}
\mbox{$\begin{picture}(22,21)(-1,0)
\put(0,10){\mb}
\put(5,10){\mb}
\put(10,10){\mb}
\put(15,10){\mb}
\put(0,-5){\mb}
\put(1,2){.}
\put(1,4){.}
\put(1,6){.}
\put(0,0){\line(0,1){10}}
\put(5,0){\line(0,1){10}}
\end{picture} $}}

\newcommand{\wyiimI}[1]
{\mbox{$\begin{picture}(20,25)(-1,0)
\put(0,0){\tiny ${k+2}$}
\put(9,7){\vector(0,1){8}}
 \put(9,-2){\vector(0,-1){8}}
\end{picture}$}
\mbox{$\begin{picture}(16,21)(-1,0)
\put(0,15){\line(1,0){10}}
\put(0,10){\line(1,0){10}}
\put(0,0){\line(1,0){10}}
\put(0,-5){\line(1,0){10}}
\put(0,-5){\line(0,1){20}}
\put(10,-5){\line(0,1){20}}
\put(5,-5){\line(0,1){7}}
\put(5,8){\line(0,1){7}}
\put(1,2){.}
\put(1,4){.}
\put(1,6){.}
\put(6,2){.}
\put(6,4){.}
\put(6,6){.}
\put(0,-10){\mb}
\end{picture}$} 
}

\newcommand{\yiII}{\mbox{$
\begin{picture}(7,10)(-1,1)
\put(0,5){\mb}
\put(0,0){\mb}
\end{picture}
$}}

\newcommand{\yiiII}{\mbox{$
\begin{picture}(13,10)(-1,1)
\put(0,10){\line(1,0){10}}
\put(0,5){\line(1,0){10}}
\put(0,0){\line(1,0){10}}
\put(0,0){\line(0,1){10}}
\put(10,0){\line(0,1){10}}
\put(5,0){\line(0,1){10}}
\end{picture}
$}}

\newcommand{\yiIIIiII}{\mbox{$
\begin{picture}(13,15)(-1,1)
\put(0,15){\line(1,0){10}}
\put(0,10){\line(1,0){10}}
\put(0,5){\line(1,0){10}}
\put(0,0){\line(1,0){5}}
\put(0,0){\line(0,1){15}}
\put(10,5){\line(0,1){10}}
\put(5,0){\line(0,1){15}}
\end{picture}
$}}

\newcommand{\yiIIIi}{\mbox{$
\begin{picture}(13,15)(-1,1)
\put(0,15){\line(1,0){10}}
\put(0,10){\line(1,0){10}}
\put(0,5){\line(1,0){5}}
\put(0,0){\line(1,0){5}}
\put(0,0){\line(0,1){15}}
\put(10,10){\line(0,1){5}}
\put(5,0){\line(0,1){15}}
\end{picture}
$}}

\newcommand{\yiIII}{\mbox{$
\begin{picture}(7,15)(-1,1)
\put(0,15){\line(1,0){5}}
\put(0,10){\line(1,0){5}}
\put(0,5){\line(1,0){5}}
\put(0,0){\line(1,0){5}}
\put(0,0){\line(0,1){15}}
\put(5,0){\line(0,1){15}}
\end{picture}
$}}

\newcommand{\yiiIII}{\mbox{$
\begin{picture}(13,15)(-1,1)
\put(0,15){\line(1,0){10}}
\put(0,10){\line(1,0){10}}
\put(0,5){\line(1,0){10}}
\put(0,0){\line(1,0){10}}
\put(0,0){\line(0,1){15}}
\put(10,0){\line(0,1){15}}
\put(5,0){\line(0,1){15}}
\end{picture}
$}}

\newcommand{\yiiiI}{\mbox{$
\begin{picture}(18,5)(-1,1)
\put(0,5){\line(1,0){15}}
\put(0,0){\line(1,0){15}}
\put(0,0){\line(0,1){5}}
\put(10,0){\line(0,1){5}}
\put(5,0){\line(0,1){5}}
\put(15,0){\line(0,1){5}}
\end{picture}
$}}

\newcommand{\yiiiiI}{\mbox{$
\begin{picture}(23,5)(-1,1)
\put(0,5){\line(1,0){20}}
\put(0,0){\line(1,0){20}}
\put(0,0){\line(0,1){5}}
\put(10,0){\line(0,1){5}}
\put(5,0){\line(0,1){5}}
\put(15,0){\line(0,1){5}}
\put(20,0){\line(0,1){5}}
\end{picture}
$}}

\newcommand{\yii}[1]{ \mbox{$\mbox{$\begin{picture}(11,21)(-1,0)
\put(1,3){\tiny ${#1}$}
\put(4,10){\vector(0,1){6}}
 \put(4,1){\vector(0,-1){6}}
\end{picture}$}
\mbox{$\begin{picture}(13,21)(-1,0)
\put(0,15){\line(1,0){10}}
\put(0,10){\line(1,0){10}}
\put(0,0){\line(1,0){10}}
\put(0,-5){\line(1,0){10}}
\put(0,-5){\line(0,1){20}}
\put(10,-5){\line(0,1){20}}
\put(5,-5){\line(0,1){7}}
\put(5,8){\line(0,1){7}}
\put(1,2){.}
\put(1,4){.}
\put(1,6){.}
\put(6,2){.}
\put(6,4){.}
\put(6,6){.}
\end{picture}$} 
$}}

\newcommand{\wyii}[1]{ \mbox{$\mbox{$\begin{picture}(20,21)(-1,0)
\put(0,3){\tiny ${#1}$}
\put(9,10){\vector(0,1){6}}
 \put(9,1){\vector(0,-1){6}}
\end{picture}$}
\mbox{$\begin{picture}(13,21)(-1,0)
\put(0,15){\line(1,0){10}}
\put(0,10){\line(1,0){10}}
\put(0,0){\line(1,0){10}}
\put(0,-5){\line(1,0){10}}
\put(0,-5){\line(0,1){20}}
\put(10,-5){\line(0,1){20}}
\put(5,-5){\line(0,1){7}}
\put(5,8){\line(0,1){7}}
\put(1,2){.}
\put(1,4){.}
\put(1,6){.}
\put(6,2){.}
\put(6,4){.}
\put(6,6){.}
\end{picture}$} 
$}}

\newenvironment{proof}{\begin{trivlist} \item[] {\em Proof}. }%
{\hfill $\Box$ \end{trivlist}}

\newcommand{\ul}[1]{\underline{#1}}

\newcommand{\super}[2]{^{\overbrace{\mbox{\tiny $#1$}}^{#2}}}
\newcommand{\sub}[2]{_{\underbrace{\mbox{\tiny $#1$}}_{#2}}}

\def\contr{\operatorname{contr}}

\def\sideremark#1{\ifvmode\leavevmode\fi\vadjust{\vbox to0pt{\vss
 \hbox to 0pt{\hskip\hsize\hskip1em
 \vbox{\hsize3cm\tiny\raggedright\pretolerance10000
 \noindent #1\hfill}\hss}\vbox to8pt{\vfil}\vss}}}%
\newcommand{\GL}{\operatorname{GL}}
\newcommand{\SL}{\operatorname{SL}}

\begin{document}
\headheight=6mm
\hoffset=-10mm
\title[Calculus for Quaternionic Structures]{Invariant 
Local Twistor Calculus for Quaternionic Structures and
Related Geometries 
}

\author{A. Rod Gover and Jan Slov\'ak}

\date{}
\thanks{The first author is an Australian Research Council 
QEII Research Fellow. The second
author supported by Australian Research Council, University of Adelaide, and
GACR grant Nr. 201/96/0310
}

\begin{abstract}
New universal invariant operators are introduced in a class of
geometries which include the quaternionic structures and their
generalisations as well as 4-dimensional conformal (spin) geometries. 
It is shown that, in a broad sense,  all invariants and invariant
operators arise from these universal operators and that they may be
used to reduce all invariants
problems to corresponding algebraic problems involving homomorphisms
between modules of certain parabolic subgroups of Lie groups. Explicit
application of the operators is illustrated by the construction of all
non-standard operators between exterior forms on a large class of the
geometries which includes the quaternionic structures.
\end{abstract}

\maketitle

\frenchspacing
\section{Introduction}

A {\em real almost Grassmannian structure} on a manifold $M$ 
(briefly a {\em real AG-structure})
is given by a fixed identification of the tangent bundle
$TM$ with the tensor product of two auxiliary vector bundles of dimensions
$p$ and $q$, together with the identification of their top
degree exterior powers. In the realm of Penrose's abstract index
notation, we shall express this by 
\begin{equation}\label{fund-ident}
{\cal E}^a=\ce_{A'}\otimes {\cal E}^A ={\cal E}^A_{A'}, \quad
\wedge^q{\cal E}^A\simeq \wedge^p{\cal E}_{A'} .
\end{equation}
Equivalently, this amounts to the reduction of the structure group
$\GL(pq, {\Bbb R})$ of the tangent bundle to its subgroup
$G_0=\operatorname{S}(\GL(p,{\Bbb R})\times \GL(q,{\Bbb R}))$. Thus
the complexified tangent bundle of a real AG-structure is equipped by
the reduction of its structure group to $G_0^{\Bbb
  C}=\operatorname{S}(\GL(p,{\Bbb C})\times \GL(q,{\Bbb C}))$.  There
is another class of geometries on $4m$-dimensional manifolds with
similar behaviour.  The geometries are defined by reductions of the
structure groups of the tangent bundles to the groups
$G_0=\operatorname{S}(\GL(p/2,{\Bbb H})\times \GL(q/2,{\Bbb
  H}))\subset \GL(pq,{\Bbb R})$ with $2\le p\le q$ even, and the
complexifications of their tangent bundles enjoy again the fundamental
identification (\ref{fund-ident}).  The most important algebraic
feature of the two types of the structures above is that, for each
pair $p,q$, their respective structure groups $G_0$ are the maximal
reductive parts of certain maximal parabolic subgroups $P$ in two
different real forms $G$ of the same complex semi-simple group
$G^{\Bbb C}=\SL(p+q,{\Bbb C})$. 
A geometry will be called an {\em AG-structure} if it has a  structure
group $G_0$ where $G_0$ is a maximal
reductive part of a parabolic $P\subset G$  such that $P=G\cap P^{\Bbb C}$
with $ G^{\Bbb C}=\SL(p+q,{\Bbb C}) $ and where $ P^{\Bbb C}$ is the maximal
parabolic in $ G^{\Bbb C}$ such that $G^{\Bbb C}/P^{\Bbb C}$ is the 
Grassmannian of complex $p$-planes in $ \C^{p+q}$ (with $2\leq p\leq
q$).
Thus the members of the list of all AG-structures are named
by such pairs $(G,P)$ and in fact the  real group $G$ is one of the
following: $G={\rm SL}(p+q,{\Bbb R})$ with $2\le p\le q$, $G={\rm
  SL}(p/2 +q/2,{\Bbb H})$ and $p,q$ are even, or $G={\rm SU}(p,p)$,
see Appendix \ref{A} for more details.  Henceforth $G$ and $P$ will
indicate such a pair and $G_0$ will be the reductive part of $P$. The identification
(\ref{fund-ident}) of the complexified tangent spaces is given for all
the AG-structures. The
complex almost Grassmannian structures on complex manifolds were
studied in \cite{BaiE} under the name `paraconformal manifolds'.
Similar objects were introduced earlier in \cite{Gon}, see also
\cite{Gin}.
 
The most well known examples of such structures are 4-dimensional 
conformal spin structures (here $G_0={\Bbb R}\cdot
\operatorname{Spin}(p,q,{\Bbb R})\subset
\operatorname{Spin}(p+1,q+1,{\Bbb R})$, $p+q=4$, and the complexification
$\operatorname{Spin}(6,{\Bbb C})\simeq \operatorname{SL(4,{\Bbb C})}$). 
We will extend the term `spinor' from that case and in all cases deem the
auxiliary bundles $\ce_{A'}$ and $ \ce^A$ to be {\em spinor bundles}.

The almost quaternionic structures on manifolds are 
classical 1st order G-structures, such that their structure group 
$G_0$ is the subgroup 
$\GL(m,{\Bbb H})\times_{{\Bbb Z}_2}\operatorname{Sp}(1)\subset \GL(4m,\Bbb
R)$,
see \cite{Sal}. We have to notice that the action of $G_0$ on ${\Bbb H}^m$
(i.e. the indicated embedding into the real general linear group)
is defined by the adjoint action of the block-diagonal matrices in
$\GL(1+m,{\Bbb H})$ on the block below the diagonal. 
The group $\tilde 
G_0=S(\GL(1,{\Bbb H})\times\GL(m,{\Bbb H}))$ is the 
universal cover of $G_0$ and the choice of the structure group
$\tilde G_0$ makes no difference locally. In
particular, the almost quaternionic structures belong to our class of
AG-structures. They are called {\em quaternionic\/} if they
admit a torsion-free connection. It was
pointed out in \cite{Sal}, and worked out in much detail in
\cite{BaiE}, \cite{Ba}, that these structures fit into a larger class
of geometries coming from the so called $|1|$-graded semi-simple Lie
algebras. This is exactly our point of view and the corresponding entry in
our list of pairs $(G,P)$ is that with $G=\SL(1+q/2,{\Bbb H})$, $q\ge 2$
even. 

Despite the very
transparent geometric differences between, for example, the real
almost Grassmannian structures and the almost quaternionic structures,
we will treat all these cases simultaneously. In the cases
corresponding to the `split
real form' $G= {\rm SL}(p+q,{\Bbb R})$ we will write $TM$ to mean the
usual tangent bundle while for the other cases $ TM$ will mean the
complexification of the tangent bundle.
Similarly for  $P$-modules, and the bundles they induce, 
we will take these to be real for the geometries of the split real
forms but complex for the geometries corresponding to the other
groups. With this understood we will suppress explicit reference to
the scalars concerned and write, for example, $\SL(m)$ for
 either the real or complex special linear group as required by context.
These conventions will enable us to use 
the same index formalism  for all these geometries and also enable us
to avoid
complexifying except where necessary.

Treating all such AG-structures simultaneously,
the main results we obtain are as follows: \\
$\bullet$ We construct a new invariant first order differential
operator that we call a twistor-D operator -- see definition
\ref{Defn}.  This operator may be viewed as an analogue, for these
structures, of the Levi-Civita
connection of Riemannian geometry.\\
$\bullet$ Via the twistor-D operator we construct curved analogues of
all the non-standard operators between exterior differential forms on a class
of AG-structures that includes all the quaternionic
geometries -- see theorem \ref{weaks}. \\
$\bullet$ We use the twistor-D operator to construct a module for an
appropriate parabolic subgroup $P$ such that all invariant
differential operators (linear and polynomial and up to any chosen
order) and invariants of AG-structures are equivalent to
$P$-homomorphisms from this module to irreducible $P$-modules. See in
particular theorem \ref{main}.  The implications of this are discussed
below.

We should also point out that considerable detail of a `calculus' to
enable manipulation and application of the twistor-D operator and its
accompanying machinery is presented. Most of this is strictly needed
to establish the results mentioned. However we have attempted to
present this in an explicit form that could be directly used by readers
as we believe that there are many potential applications of these tools in
mathematical-physics, especially since they include new tools for the
4-dimensional conformal structures and their associated twistor theory.
For example the twistor-D operator should be particularly useful for the
construction and study of conformally invariant spinor equations in
4-dimensions.  In addition to the main results there are other
observations and results along the way. In particular, in section
\ref{observations} we observe an obvious extension to Salamon's
complex, we relate the twistor-D operator to the so called tractor-D
operator of conformal geometry and we also generalise the latter to a class
of AG-structures.

Underlying our constructions here is the result that a manifold with
an almost Grassmannian structure comes equipped with a canonical {\em
  Cartan bundle} ${\cal G}\to M$ and associated canonical Cartan
connection. In each case $\cg$ is a principal fibre bundle with
structure group $P$ where this is a maximal parabolic subgroup of a
Lie group $G$ as above. The canonical {\em (normal) Cartan connection}
$\omega$ is a special 1-form on $\cg$ which takes values in the Lie
algebra $\frak g$ of $ G$ and gives a complete parallelisation of $\cg
$ (see Appendix \ref{A}, \cite{Ba} and \cite{CSS2}, for more deails).
The Cartan bundle may be regarded as a deformation of the homogeneous
situation where one has $G$ as principal bundle with fibre $ P$ over
$G/P$ and in this latter picture the Cartan connection reduces to the
Maurer Cartan form.  As in the homogeneous case each $P$-module ${V}$
gives rise to an induced or {\em natural bundle} ${\cal V}$.
Moreover, in the special case of a $G$-module $W$ the corresponding
natural bundle ${\cal W}$ comes equipped with a canonical linear
connection (also denoted by $\omega$). Such bundles will be described 
here by what will be 
called {\em (local) twistor bundles} and their canonical connections
will be viewed as {\em twistor connections}.  Since, in the current 
work, we are
concerned with the production of explicit operators on $M$ we avoid
a detailed discussion of the Cartan bundle and work directly on these
induced natural bundles and their connections. Indeed most of the work
can be understood without a deep understanding of the inducing
Cartan bundles. However we would like to point out that many of the
`background results' can be recovered most efficiently from the 
principal bundle
point of view and Appendix \ref{A} is dedicated to extracting from the
general theory of Cartan bundles and their connections (as in for
example \cite{CSS1}) the results required for the current work.

Calculus similar to the twistor calculus we develop here has been
successfully applied to other related geometries. For example in
\cite{Goproj} a first order invariant tractor-D operator (rediscovered
in \cite{BaiEGo} but originally due to Tracey Thomas), and some calculus based
around this, is used to construct all density valued invariants of
projective geometries. In \cite{Goconf} a similar programme is in place
to produce a complete invariant theory for conformal geometries and
there are already many new results in this.  Such calculus has also
been used to proliferate invariant operators on conformal, projective
and CR structures. As with the AG-structures studied here, these geometries are
all `parabolic geometries' which may be viewed as deformations of
homogeneous structures $G/P$ where $ G$ is semi-simple and $ P$ a
parabolic subgroup. It turns out that at each point of such a
structure $ P$ acts on the jet information (jets of the geometric
structure itself or jets of a field on the structure). Understanding
and dealing with this action is the key problem. This is difficult and
subtle in general and many papers have discussed similar problems, see
e.g. \cite{Ba}, \cite{CSS1}, \cite{Sl}, and the references therein.
Roughly speaking the programme here, as with the tractor calculus, is
to use the twistor-D operator to package this jet information into
`parcels' which are $P$-submodules of irreducible $G$-modules. This is
a huge step since at least the latter $G$-modules are understood and
can be dealt with by classical techniques such as Weyl's invariant
theory. (A discussion of  the general programme, in the context of
tractor calculus, as well as other results are described in  
\cite{Gosrni}.) Then invariants and invariant operators may be
proliferated by identifying 
the relevant $P$-submodules of irreducible $G$-modules.
That all
invariants and invariant operators are equivalent to the corresponding
$P$-homorphisms is the content of theorem \ref{main}. A more intuitive
interpretation of this result is that all invariants arise from the
twistor-D operator (and its concatenations -- the universal invariant
$\Do^{(k)}$ operators of section \ref{twistcalc}).  As far as we know
this is the first theorem of its sort and thus far there is no
corresponding theorem established for the tractor operators.  This
theorem leaves open the question of whether the remaining
$P$-submodule problems are tractable. Evidence that in many important
cases they are is the success of the analogous tractor calculus, as
mentioned above, and more importantly for this case the application of
the twistor-D operator to produce the new family of invariant
differential operators in section \ref{newops}. For future work in
this direction, as well as to develop some results needed here,
appendix \ref{compseries} discusses the composition series of submodules
in a rather general setting.

The plan of the paper goes as follows. After setting notation and
outlining some further preliminaries in the next section, we introduce
the twistor-D operators.  The fourth section is devoted to the main
Theorem \ref{main} the proof of which relies on an explicit
description of the normal forms of the AG-structures, cf. Appendix
\ref{normal}. Then we proceed with our main application, the curved
analogues of the non-standard operators on 
exterior forms. These are fourth order and include analogies to the
square of the Laplacian in 
four-dimensional conformal geometries. Further observations, as
mentioned above, in  
Section \ref{observations}, are followed by the three Appendices.

\smallskip
\noindent{\bf Acknowledgements.} Discussions with Andreas \v Cap and
Michael Eastwood were important. Experimenting with
Brian Boe's computer program for computing the classification 
lists of homomorphisms between generalized Verma modules has been also
very useful (cf. \cite{BoeC}). Essential parts of the research were 
done during the second author's stays at the University of Adelaide and QUT
in Brisbane, and the first author's stay at Masaryk University in Brno. Some
writing was also done during the authors' visit at Erwin Schr\"odinger
Institute in Vienna. 

\section{Preliminaries}\label{preliminaries}

Here we review some important technicalities and introduce our
notational conventions. We omit explicit verifications of most of the 
claims as they follow easily from the general theory as reviewed in
Appendix \ref{A}, see also  \cite{CSS1}. For an explicit development
(although in the complex setting)
with notation and conventions very similar to those here see \cite{BaiE}.

\smallskip
\noindent {\bf Index formalism.}
Except where otherwise indicated we use Penrose's abstract index
notation \cite{ot} which allows for easy explicit calculations without
involving a choice of basis.  Thus we may write, $v^A$ or $ v^B$ for a
section of the unprimed fundamental spinor bundle $\ce^A$. Similarly
$w_{A'}$ could denote a section of the primed fundamental spinor
bundle $\ce_{A'}$.  We write $\ce_A$ for the dual bundle to $\ce^A$
and $\ce^{A'}$ for the dual to $\ce_{A'}$.  The tensor products of
these bundles yield the general spinor objects such as
$\ce_{AB}:=\ce_A\otimes\ce_B$, $\ce^{ABC'}_{A'B'}$ and so forth.  The
tensorial indices are also abstract indices. Recall (see above)
that $\ce^a=\ce^A_{A'}$ is the tangent bundle, so
$\ce_a=\ce^{A'}_A$ is the cotangent bundle and we may use the terms
`spinor' or `section of a spinor bundle' to describe tensor fields. 
  
  
A spinor object on which some indices have been contracted will be
termed a {\em contraction} (of the underlying spinor). For example
$$
v^{ABC'}_{BC'DE} 
$$ 
is a contraction of $v^{ABC'}_{DD'EF}$. In many cases the underlying
spinor of interest is a tensor product of lower valence spinors. For
example
$$
v^{AB}w^{C'}_Bu_{ACD}
$$
is a contraction of $v^{AB}w^{C'}_C u_{DEF}$. The same conventions are
used for the tensor indices and the twistor indices; the latter are to
introduced below. 
Standard notation is also used
for the symmetrizations and antisymmetrizations over some indices.

\smallskip
\noindent{\bf Weights and scales.} We define line bundles of densities or
{\em weighted functions} as follows.
The weight $-1$ line bundle $\ce[-1]$ over $M$ is identified with 
$$
\ce\super{[A'B'\cdots C']}{p}.
$$ 
Then, for integral $w$, the weight $w$ line bundle $\ce[w]$ is defined to be
$(\ce[-1])^{-w}$. In fact in the case of AG-geometries corresponding to
the real-split form SL$(p+q,{\Bbb R})$ we can (locally) extend this definition
to weights $w\in  \br$ by locally selecting a ray fibre subbundle of
$\ce[-1]$. Calling this say $\ce_+[-1]$ we can then define the ray
bundles  $\ce_+[w]: =(\ce[-1])^{-w}$. Finally these may be canonically
extended to line bundles in the obvious way. 
In any case we write $\ce^{A'}[w]$ for $\ce^{A'}\otimes \ce[w]$ and
so on, whenever defined. 
In view of the defining isomorphism 
\begin{equation}\label{defisom}
h: \wedge ^q\ce^A \toiso \wedge^p\ce_{A'}
\end{equation} 
we also have
$$
\ce[-1]\cong\ce\sub{[AB\cdots C]}{q},\quad 
\ce[1]\cong\ce\super{[AB\cdots C]}{q}\cong\ce\sub{[A'\cdots C']}{p}
.$$
We write $\vol^{A'B'\cdots C'}$ for the tautological section of 
$\ce^{[A'B'\cdots C']}[1]$ giving the mapping $\ce[-1]\toiso \wedge^p\ce^{A'}$
by
\begin{equation}\label{raiselower}
f\mapsto f\vol^{A'B'\cdots C'},
\end{equation}
and $\vol_{D\cdots E}$ for similar object giving $\ce[-1]\toiso
\wedge^q \ce_A$.
A {\em scale} for the AG-structure is a nowhere 
vanishing section $\xi$ of
$\ce[1]$.  Note that such a choice is equivalent to a choice of spinor 
`volume' form
$$
\vol_\xi^{A'\cdots C'}:=\xi^{-1}\vol^{A'\cdots C'},
$$
{\em or} to a choice of form,
$$
\vol^\xi_{D\cdots  E}:=\xi^{-1}\vol_{D\cdots E} .
$$

\smallskip
\noindent{\bf Distinguished connections.}
A connection $\nd_a$ on $M$ belongs to
the given AG-structure (this really means $\nd_a$ comes from a principal 
connection on the bundle $\cg_0$ described below) 
if and only if it satisfies two conditions: 
\begin{itemize}
\item $\nd_a$ is the tensor product of linear connections (both of
  which we shall also denote $\nd_a$)
on the spinor bundles $\ce^A$ and $\ce_{A'}$, 
\item the defining isomorphism $h$ in (\ref{defisom}) is covariantly constant, 
i.e. $\nd_a h=0$.
\end{itemize}
Our conventions for the torsion $T_{ab}{}^c$  and curvature
$R_{ab}^{\vphantom{b}}{}^c_d$ of a connection $\nd_a$ on
the tangent bundle $TM$ are determined by the following equation,
$$
2\nd_{[a}\nd_{b]}v^c=T_{ab}{}^d\nd_d v^c + R_{ab}^\vb{}^c_d v^d .
$$

Since  $T_{ab}{}^c$ is skew on its lower indices,
$T_{ab}{}^c= T_{[ab]}{}^c$, it can be written as a sum of two 
terms
$$
T_{ab}{}^c=F_{ab}{}^c+\tilde{F}_{ab}{}^c
$$
where 
$$
F_{ab}{}^c=:F^{A'B'C}_{ABC'}=F_{(AB)C'}^{[A'B']C}~\mbox{ and }~ 
\tilde{F}_{ab}{}^c=: \tilde{F}^{A'B'C}_{ABC'}=F_{[AB]C'}^{(A'B')C}
.$$

The Cartan bundle $\cg$ over the manifold $M$ has the
quotient $\cg_0$, a principal fibre bundle with structure group $G_0$. By
the general theory, each $G_0$-equivariant section $\sigma:\cg_0\to \cg$ of
the quotient projection defines the distinguished principal connection on 
$\cg_0$, the pullback of the ${\frak g}_0$-part of $\omega$. The whole class of
these connections consists precisely of connections on $\cg_0$ with the unique 
torsion taking values in the kernel of $\partial^*$. A straightforward 
computation shows that the latter condition is equivalent to 
the condition that both $\tilde{F}$
and $F$ be completely trace-free (cf. Appendix \ref{A}). 
Each principal connection on $\cg_0$ induces the induced
connection on the bundle $\ce[1]\setminus \{0\}$ which is associated to
$\cg_0$ and, moreover, the resulting correspondence between the sections
$\sigma$ and the latter connections is bijective.  
In particular, each section $\xi$  
of the bundle $\ce[1]\setminus\{0\}$ defines
uniquely a reduction $\sigma$, such that the corresponding distinguished
connection leaves $\xi$ horizontal. Altogether we have
recovered Theorems 2.2, 2.4 of \cite{BaiE}. We rephrase these here for
convenience: 
\begin{theorem} 
Given a scale $\xi$ on an AG-structure there are unique connections on
$\ce^A$ and $\ce_{A'}$ such that $F_{ABC'}^{A'B'C}$  and 
$\tilde{F}_{ABC'}^{A'B'C}$ are totally trace-free, the
induced covariant derivative preserves the isomorphism $h$ of
\nn{defisom}, and $\nd_a \xi=0$. The torsion components $F_{ab}{}^c$ and
$\tilde{F}_{ab}{}^c$ of the induced connection on $TM$ 
are invariants of the AG-structures.
\end{theorem}

Notice that in the special case of the four-dimensional conformal
geometries, there is always a connection with vanishing torsion on
$\cg_0$ and so both $F$ and $\tilde F$ are zero. The scales correspond
to a choice of metric from the conformal class while the general
distinguished connections (corresponding to the reduction parameter
$\sigma$ being not necessarily
exact) are just the Weyl geometries.
  
We may write $\nd^{\xi}_a$ to
indicate a connection as determined by the theorem, although mostly
we will omit the $\xi$. Thus
we might write $\nd^{\hat\xi}_a$ or simply $\hat{\nd}_a$ to indicate
the connection corresponding to a scale $\hat\xi$ and similar
conventions will be used for other operators and tensors that depend on $\xi$. 

In what follows, for the purpose of explicit calculations, we shall
often 
choose a scale and work with the corresponding connections.
Objects are then well defined, or {\em invariant}, on the AG-structure
if they are independent of the choice of scale.
Note that if we change the scale according to $\xi\mapsto
\hat{\xi}=\Omega^{-1}\xi$, where $\Omega$ is a smooth non-vanishing
function, then the connection transforms as follows:
\begin{equation}\label{ndtrans}
\begin{array}{lcl}
\ce^A &:& \hat\nd^{A'}_A u^C= \nd^{A'}_A u^C+\delta^C_A\up^{A'}_B u^B \\
\ce_{A'} &:&\hat\nd^{A'}_A u_{C'}= \nd^{A'}_A u_{C'} +\delta^{A'}_{C'}\up^{B'}_Au_{B'}\\
\ce_B &:&\hat\nd^{A'}_A v_B=\nd^{A'}_A v_B -\up^{A'}_Bv_A\\ 
\ce^{B'} &:&\hat\nd^{A'}_A v^{B'}= \nd^{A'}_A v^{B'}- \up^{B'}_A v^{A'}\\
\end{array}
\end{equation}
where $\up_a:=\Omega^{-1}\nd_a \Omega$.
Consequently 
\begin{equation}\label{wttr}
\widehat{\nd}_a f=\nd_a f+ w \up_a f 
\end{equation}
if $f\in \ce[w]$. All these formulae follow from the general discussion in
Appendix \ref{A}, but they are also easily checked directly.

Given a choice of scale $\xi$, we will write $R_{ab}^\vb{}^C_D$ (or
$R^{(\xi)}_{ab}{}^C_D$ to emphasise the choice of scale)
for the curvature of $\nd_a$ on $\ce^A$ and $R_{ab}^\vb{}^{C'}_{D'}$ for
the curvature of $\nd_a$ on $\ce_{A'}$, that is
$$
(2\nd_{[a}\nd_{b]}-T_{ab}{}^e\nd_e)v^C=R_{ab}^\vb{}^C_D v^D,\ 
(2\nd_{[a}\nd_{b]}-T_{ab}{}^e\nd_e)w_{D'}=-R_{ab}^\vb{}^{C'}_{D'}w_{C'}.
$$
Then the curvature of the induced linear connection on $TM$ is  
$$
R_{ab}^\vb{}^c_d=R_{ab}^\vb{}^{C'}_{D'}\delta^C_D+R_{ab}^\vb{}^C_D\delta^{C'}_{D'}
.$$
Observe that since $\nd_a$ preserves the volume forms 
$\vol^\xi_{A'\cdots C'}$ and $\vol^\xi_{D\cdots E}$ it
follows that $R_{ab}^\vb{}^C_D$ and 
$R_{ab}^\vb{}^{C'}_{D'}$ are trace-free on the spinor indices
displayed. Thus the equations
$$
R_{ab}^\vb{}^C_D=
U_{ab}^\vb{}^C_D-\delta^C_B\Rho^{A'B'}_{AD}+\delta^C_A\Rho^{B'A'}_{BD}
$$
and 
$$
R_{ab}^\vb{}^{C'}_{D'}=
U_{ab}^\vb{}^{C'}_{D'}+\delta^{B'}_{D'}\Rho^{A'C'}_{AB}-
\delta^{A'}_{D'}\Rho^{B'C'}_{BA}
$$
determine the objects $U_{ab}^\vb{}^C_D$, $U_{ab}^\vb{}^{C'}_{D'}$ and the
{\em Rho-tensor},
$\Rho_{ab}$, if we require that $U_{ACD}^{A'B'C}=0=U_{ABD'}^{A'D'C'}$.  
In this notation we have,
\begin{equation}\label{curvdec}
R_{ab}^\vb{}^c_d= U_{ab}^\vb{}^c_d
+\delta^{D'}_{C'}\delta^C_A\Rho^{B'A'}_{BD}-
\delta^{D'}_{C'}\delta^C_B\Rho^{A'B'}_{AD}-
\delta^{C}_{D}\delta^{A'}_{C'}\Rho^{B'D'}_{BA}
+\delta^{C}_{D}\delta^{B'}_{C'}\Rho^{A'D'}_{AB} 
\end{equation}
where
\begin{equation}\label{Wsplit}
U_{ab}^\vb{}^c_d=U_{ab}^\vb{}^C_D\delta^{D'}_{C'}+U_{ab}^\vb{}^{D'}_{C'}\delta^C_D
.
\end{equation}
In the case of $p=2=q$ this agrees with the usual
decomposition of the curvature of the Levi-Civita connection into 
the conformally invariant (and trace-free) Weyl tensor part and the
remaining part given by the Rho-tensor (see e.g. \cite{BaiEGo}). 
All these equations also follow from the general definitions of the $U$'s
and $\Rho$'s in (\ref{A-U}). Note that $U$'s are two-forms
valued in ${\frak g}_0$ coming from the curvature of the canonical Cartan
connection and so they are in the kernel of $\partial ^*$. This is the
source of the condition on the trace, but they are not trace-free in
general: 
\begin{equation}\label{skewrho1}
U_{ab}^\vb{}^C_C=-U_{ab}^\vb{}^{C'}_{C'}=2\Rho_{[ab]}
.\end{equation}
On the other hand, it follows from the Bianchi identity,
$$
R_{[ab}^{}{}_{c]}^d+\nd_{[a}T_{bc]}{}^d +T_{[ab}{}^eT_{c]e}{}^d=0,
$$
that 
\begin{equation}\label{skewrho}
2(p+q)\Rho_{[ab]}=-\nd_cT_{ab}{}^c .
\end{equation}
The Rho-tensor $\Rho_{ab}$ has the transformation equation
\begin{equation} \label{Ptrans}
\hat{\Rho}{}_{AB}^{A'B'}=\Rho_{AB}^{A'B'}-\nd_{A}^{A'}\Upsilon_{B}^{B'}
+\Upsilon_{A}^{B'}\Upsilon_{B}^{A'} .
\end{equation}
Again, this can be easily checked directly but we give a general explanation
in (\ref{rho}). 

We are most interested in the special case $p=2$. Then the 
whole component $F_{ab}{}^c$ is irreducible and
so it vanishes by our condition on the trace, while the other
component $\tilde{F}_{ab}{}^c$ of the torsion, together with the trace-free
part of $U_{(ABC)}^{[A'B']D}$ are the only local invariants of 
the structures (i.e. the
AG-structure is locally flat if and only if these two vanish). In all other
cases $2<p\le q$, the two components of the torsion are the only invariants,
cf. the end of Appendix \ref{A}.
 
The totally symmetrized covariant derivatives of the Rho-tensors will play a
special role. We will use the notation 
$$
S_{a\cdots b}:=\underbrace{\nd_{(a}\nd_b\cdots \nd_d \Rho_{ef)}}_s
$$
for $s=2,3\cdots$.

\smallskip
\noindent{\bf Twistors.}
Via the Cartan bundle $\cg$ over $M$ 
any $P$-module $V$
gives rise to a {\em natural bundle} (or induced bundle)  
$\cv$. Sections of $\cv$ are identified with 
functions $ f:\cg \to V $ such that $f(x.p)=\rho(p^{-1})f(x)$, 
where $x\mapsto x.p$ gives the  action of $p\in P $ on $ x\in \cg$ while
$\rho$ is the action defining the $P$-module structure.  

Recall also that the Cartan bundle is equipped with a canonical
connection, the so called normal Cartan connection $\omega$. In view
of this it is in our interests to work, where possible, with natural
bundles $\cv$ induced from $V$ where this is not merely a $ P$-module
but in fact a $G$-module. Then the Cartan connection induces an
invariant linear connection on $\cv$.  Let us write $V^\alpha$ for the
module corresponding to the standard representation of
$G$ on ${\Bbb R}^{p+q}$ and write $V_\alpha$ for the dual module.  The
index $\alpha$ is another Penrose-type abstract index and we write
$\ce^\alpha$ and $ \ce_\alpha$ for the respective bundles induced by
these $G$-modules.  All finite dimensional $G$-modules are submodules
in tensor products of the fundamental representations $V^\alpha$ and
$V_\alpha$.  Thus the bundles $\ce^\alpha$ and $ \ce_\alpha$ play a
special role and we term these {\em (local) twistor bundles} (c.f. \cite{BaiE,nt}). In fact in line with the use of the word
``tensor'' we will also describe any explicit subbundle of a tensor
product of these bundles as a twistor bundle and sections of such
bundles as local twistors.
In particular observe that there is a canonical completely skew
local-twistor $(p+q)$-form 
$h_{\alpha\beta\cdots \gamma}$
on $\ce^\alpha$ which is equivalent to the isomorphism \nn{defisom}.
We write $h^{\alpha\beta\cdots \gamma}$ for the dual
completely skew twistor satisfying
$h^{\alpha\beta\cdots \gamma}h_{\alpha\beta\cdots \gamma}=(p+q)!$.

All finite dimensional $P$-modules enjoy filtrations which split
completely as $G_0$-modules.  $ V^\alpha$ and $ V_\alpha$, give the
simplest cases and, as $P$-modules, admit filtrations
$$
V^\alpha = V^A+V^{A'},\quad V_\alpha=V_{A'}+V_A .
$$
(Our
notational convention is that the `right ends' in the formal sums are
submodules while the `left ends' are quotients.) These determine 
filtrations of the twistor bundles 
$$
\ce^\alpha =\ce^A + \ce^{A'}, \quad \ce_\alpha =\ce_{A'} + \ce_A .
$$
We write $X^\alpha_{A'}$ for the canonical section of
$\ce^\alpha_{A'}$ which gives the injecting morphism $\ce^{A'}\to
\ce^\alpha$ via
\begin{equation}\label{st1}
v^{A'}\mapsto X^\alpha_{A'}v^{A'}.  
\end{equation}
Similarly $Y_\alpha^A$ describes the injection of $\ce_A$ into dual
twistors,
\begin{equation}\label{st2}
\ce_A\ni u_A\mapsto Y^A_\alpha u_A \in \ce_\alpha .
\end{equation}

It follows from standard representation theory that a choice of
splitting of the exact sequence,
$$
0\to V^{A'} \to V^\alpha \to V^A\to 0
$$
is equivalent to the choice of subgroup of $P$ which is isomorphic to 
$G_0$. It follows immediately that a choice of splitting of the
twistor  bundle $\ce^\alpha$ is equivalent to a reduction from $ \cg$
to $ \cg_0$. Such a splitting  is a $G_0$-equivariant homomorphism 
$\xi:\ce^\alpha\to \ce^{A'}$. We can regard $\xi$ here as a section of
$\ce_\alpha\otimes\ce^{A'}=\ce_\alpha^{A'}$ and then in our index
notation the homomorphism is determined  by $v^\alpha\mapsto
\xi^{A'}_\alpha v^\alpha$, for any section $ v^\alpha$ of $
\ce^\alpha$. The composition of $ \xi$ with the monomorphism $
\ce^{A'}\to \ce^\alpha$ must be the identity so we have,
$$
\xi^{A'}_\beta X^\beta_{B'}=\delta^{A'}_{B'}.
$$
A splitting $\xi^{A'}_\alpha$ of $\ce^\alpha$ determines a dual
splitting $\lambda^\alpha_{A}$ of $\ce_\alpha$,
$\lambda^\alpha_{A}:\ce_\alpha\to \ce_A$. Given such splittings we
have $\ce_\alpha =\ce^A \oplus\ce^{A'}$ and 
$\ce_\alpha =\ce_{A'} \oplus \ce_A  $, so  we may write sections of these
bundles as a ``matrices'' such as
$$
[u^\alpha]_\xi=\left(\begin{array}{c}u^A\\
                               u^{A'}\end{array}\right)
\in [\Gamma\ce^\alpha]_\xi \quad
[v_\alpha]_\xi=(v_A \:\: v_{A'})\in[\Gamma\ce_\alpha]_\xi. 
$$
We will always work with splittings determined by a choice of scale
$\xi\in\ce[1]$, as discussed earlier. (That we have used the same
symbol as used for the kernel part of the symbol for the splitting is 
of course no
accident. In fact the direct connection between the scale $\xi$ and
the corresponding section $\xi^{A'}_\alpha$ is given explicitly on page
\pageref{xistuff}.) 
If $u^\alpha$ and $ v_\alpha$, as displayed, are
expressed by such a scale then the change of scale $\xi\mapsto
\hat{\xi} = \Omega^{-1}\xi$ yields a transformation of these
splittings. For example $[u^\alpha]\mapsto [u^\alpha]_{\widehat{\xi}}$
where
$$
[u^\alpha]_{\widehat{\xi}}=\left(\begin{array}{c}\hat{u}^A\\
                               \hat{u}^{A'}\end{array}\right)
              =\left(\begin{array}{c}u^A\\
                               u^{A'}-\up^{A'}_{B}u^B \end{array}\right). 
$$
With this understood we will henceforth drop the notation
$[\cdot]_\xi$ and simply write, for example, $v_\alpha\mapsto
\hat{v}_\alpha$ where
$$
\hat{v}_\alpha=(\hat{v}_A \:\: \hat{v}_{A'})=(v_A+\up^{B'}_Av_{B'} \:\: v_{A'}),
$$
for the corresponding transformation of $v^\alpha$.
In particular, the objects $\xi^{B'}_\alpha$, $\lambda^\beta_A$
are not invariant and
$$
\hat{\xi}^{B'}_\alpha =\xi_\alpha^{B'}-Y^A_\alpha\Upsilon^{B'}_A
\hspace{1.5cm}
\hat{\lambda}^\beta_A=\lambda^\beta_A + X^\beta_{B'}\Upsilon^{B'}_A .
$$
However, note that, in the splittings they determine,
$\xi^{B'}_\alpha$ and $\lambda^\beta_A$
are given
$$
\xi_\alpha^{B'}=\left( 0 \:\: \delta^{B'}_{A'} \right)
\quad
\lambda^\beta_A=\left(\begin{array}{c}
                                      \delta^B_A \\
                                         0        \end{array}\right) .
                                     $$
                                     In any such splitting the
                                     invariant objects $X^\alpha_{B'}$
                                     and $Y^A_\beta$ are given by
$$
X^\alpha_{B'} =\left(\begin{array}{c}    0\\ 
                                    \delta^{A'}_{B'}    \end{array}\right)
\quad       Y^A_\beta=\left(\delta^{A}_{B}\:\: 0 \right) .
$$ 
The first four identities of the following display are immediate,
while the final two items are useful definitions:
\begin{equation}\label{algids}
\begin{array}{lclclcl}
Y^A_\beta X^\beta_{A'}&=&0& & \xi^{A'}_\beta\lambda^\beta_A&=&0\\
Y^A_\beta \lambda^\beta_{B}&=&\delta^A_B&
                         & \xi^{A'}_\beta X^\beta_{B'}&=&\delta^{A'}_{B'}\\
Y^A_\beta \lambda^\gamma_{A}&= :&\lambda^\gamma_\beta&
                         & \xi^{A'}_\beta X^\gamma_{A'}&= :&\xi^{\gamma}_{\beta}
\end{array}
\end{equation}


We shall mostly deal with {\em weighted twistors}, i.e. tensor
products of the form $\ce^{\alpha\dots\beta}_{\gamma\dots\delta}[w]=
\ce^{\alpha\dots\beta}_{\gamma\dots\delta}\otimes\ce[w]$. 
All the above algebraic machinery
works for the weighted twistors.
In fact we
shall often omit the word `weighted' even though, of course, these
bundles do not come from $G$-modules for $w\ne0$.

Finally, we observe that via this machinery any spinorial quantity may be
identified with a (weighted) twistor. 
For example valence 1 spinors in $\ce^{A'}[w_1]$ or $\ce_{A}[w_2]$ may
be dealt
with via \nn{st1} or \nn{st2} respectively. This determines an
identification for tensor powers by treating each factor in this
way. This does all cases since, via \nn{raiselower}, 
$$
\ce^A\cong \ce\sub{[B\cdots D]}{q-1}[1] ,\quad\ce_{A'}\cong
\ce\super{[B'\cdots C']}{p-1}[-1] .$$
Now, any irreducible
representation of $G_0$ is given as a tensor product of two
irreducible components in tensor products of the fundamental spinors
(viewed as representations of the special linear groups, adjusted by a
weight).  Applying the corresponding Young symmetrizers \cite{ot,FH}
to the tensor products of $\ce_\alpha$ and $\ce^\beta$, we obtain the
explicit realization of each irreducible spinor bundle as the
subbundle of the (weighted) twistor bundle which is isomorphic to the
injecting part (in the composition series -- see appendix
\ref{compseries}) of the twistor bundle.  
Thus a section of a weighted irreducible spinor bundle $\cv$
 may be identified with a twistor object which is zero in all its
 composition factors except the first. So, in fact,  this non-zero factor is
  also the projecting part of the twistor. We write
 $\tilde{\cv}$ for this twistor (sub-)bundle satisfying $\cv\cong
 \tilde{\cv}$. Altogether, we
have established the following result.
\begin{lemma}\label{spintw}
Any irreducible spinor object $v$ can be identified with the twistor
$\tilde{v}$ which has the spinor as its projecting part. This identification
is provided in a canonical algebraic way.
\end{lemma}
In this connection we may also talk about the algebraic construction
providing the twistor bundle $\tilde{\cv}$. In any concrete case the
identifications may be described explicitly and in a rather obvious
way using 
the projectors $X, Y, \lambda, \xi$.

\section{Twistor calculus}\label{twistcalc}

Given a choice of scale $\xi$, a {\em twistor connection} $\nd_a$ 
on $\ce^\alpha$ and
$\ce_\alpha$ is given by the following formulae: 
\begin{equation}\label{twconn}
\nd^{P'}_A \left(\begin{array}{c}v^B\\
                                 v^{B'} \end{array}\right)=
\left(\begin{array}{c} \nd^{P'}_A v^B +\delta^B_A v^{P'}\\
              \nd^{P'}_A v^{B'}-\Rho^{P'B'}_{AB}v^B \end{array}\right)
\end{equation}
and
\begin{equation}\label{twconn2}
\nd^{P'}_A (u_B \:\: u_{B'})=(\nd^{P'}_A
u_B+\Rho^{P'B'}_{AB}u_{B'}\:\:\:
\nd^{P'}_A u_{B'}-\delta^{P'}_{B'}u_A ),
\end{equation}
(c.f. \cite{PenMac,Dighton,BaiE}).
Notice that whereas on the left hand side $\nd$ indicates the twistor
connection, on the right hand side the symbol $\nd$ indicates the
usual spinor connection determined by the choice of scale.
Although we have fixed a choice of scale to present explicit formulae
for these connections,
it is easily verified directly using the formulae \nn{ndtrans} that
the twistor connections  are in fact independent of the choice of
scale and so are invariant operators on the AG-structure.

An easy calculation reveals that 
\begin{multline*}
([\nd_a ,\nd_b]-T_{ab}{}^d\nd_d) \left(\begin{array}{c}v^C\\
                                 v^{C'} \end{array}\right)=
\\=
\left(\begin{array}{c} U_{ab}^\vb{}^C_D v^D -T_{ab}{}^C_{D'}v^{D'} \\[2mm]
      -2\nd_{[a}\Rho_{b]}{}^{C'}_{D}v^D
                                 +T_{ab}{}^E_{E'}\Rho^{E'C'}_{ED}v^D
                                 +U_{ab}^\vb{}^{C'}_{D'}v^{D'}
                                 \end{array}\right) .
\end{multline*}
Thus the curvature of the twistor connection is given, in this scale,
by
\begin{equation}\label{twcurv}
W_{ab}^\vb{}^\gamma_\delta=\left(\begin{array}{cc} U_{ab}^\vb{}^C_D 
& -T_{ab}{}^C_{D'} \\
-2Q_{ab}{}^{C'}_{D}    
& U_{ab}^\vb{}^{C'}_{D'}
    \end{array}\right), 
\end{equation}
where
$$
Q_{abc}:=\nd_{[a}\Rho_{b]c}-\frac{1}{2}T_{ab}{}^e\Rho_{ec} . 
$$
Note that since the twistor connection is invariant it follows that
this {\em twistor curvature} $W_{ab}^\vb{}^\gamma_\delta$ is
invariant. In fact, viewed as a ${\frak
g}$-valued 2-form on the Cartan bundle $\cg$, this is just the
curvature of the normal Cartan connection. In particular, we know that the
structures are torsion-free (in the sense of the Cartan connection) 
if and only if the torsion part 
$T_{ab}{}^c$
vanishes and they
are locally flat if and only if the whole $W_{ab}^\vb{}_{\delta}^\gamma$
vanishes.

%

\smallskip
\noindent{\bf The D-operators.} 
Observe that if $f\in\ce [w]$ then it follows easily from \nn{wttr}
that the spinor-twistor object
$$
\mbox{$D^{A'}_\beta f:= (\nd^{A'}_B f \:\: w\delta^{A'}_{B'}f)$}
$$
is invariant. 
We may regard this as an
injecting part 
of the invariant twistor object $D^\alpha_\beta
f:=X^\alpha_{A'}D^{A'}_\beta f$. By regarding, in this formula for
$D^\alpha_\beta$, $\nd$ to be the coupled twistor-spinor connection it
is easily verified that the operator $D^\alpha_\beta$ is well defined
and invariant on sections of the weighted twistor bundles
$\ce_{\alpha\cdots \gamma}^{\rho\cdots \mu}[w]$.

\begin{definition} \label{Defn}
The invariant operators $D^\alpha_\beta:
\ce_{\delta\cdots \gamma}^{\rho\cdots \mu}[w]\to 
\ce_{\beta\delta\cdots \gamma}^{\alpha\rho\cdots \mu}[w]$ are called the  
\underline{twistor-D} operators.
\end{definition}
For many calculations, where a choice of scale is made, 
it is useful to allow $D^\alpha_\beta$
to operate on spinors and their tensor products, although in this case
the result is not independent of the scale.  For example, if 
$v_C\in\ce_C[w]$ then 
$$
\mbox{$D^{A'}_\beta v_C:= (\nd^{A'}_B v_C \:\: w\delta^{A'}_{B'} v_C)$}..
$$

Since the operator
$D^\alpha_\beta$ and its concatenations will have an important role 
in the following
discussions we develop notation for their target spaces. First let
$\cf^\rho$ be defined as follows,
$$
\cf^\rho:=\mbox{ker}(Y^A_\rho: \: \ce^\rho\to\ce^A).
$$
Then we write 
$$
\cf^{\rho\cdots \sigma}_{\alpha\cdots
  \beta}:=\cf^\rho\otimes\cdots\otimes\cf^\sigma\otimes\ce_\alpha\otimes\cdots
\otimes\ce_\beta ,
$$
and $\cf^{\rho\cdots \sigma}_{\alpha\cdots\beta}[w]=
\cf^{\rho\cdots
  \sigma}_{\alpha\cdots\beta}\otimes \ce[w]$. Finally let
$$
\cs^{\mbox{\tiny $\rho\cdots \sigma$}}\sub{\alpha\cdots\beta}{k}[w]:=(\odot^k
\cf^\rho_\alpha)\otimes \ce[w].$$ 
 Note that sections of 
$\cf^\alpha_\rho$ ($=\cs^\alpha_\rho$) are not 
generally trace-free, but that
$\cf^\alpha_\rho$ is in a complement to  the trace-part of
$\ce^\alpha_\rho$.

Now if $f\in \ce[w]$ then $D^\rho_\alpha f\in \cf^\rho_\alpha
[w]$. Similarly observe  
that if $v^\sigma\in\cf^{\sigma}$ 
then 
$$
D^\rho_\alpha v^\sigma -\delta^\sigma_\alpha v^\rho
$$
is in $\cf^{\rho\sigma}_\alpha$. Thus 
$$
\Do^{\rho\sigma}_{\alpha\beta} : = \frac{1}{2}(D^\rho_\alpha D^\sigma_\beta 
+D^\sigma_\beta D^\rho_\alpha - \delta^\sigma_\alpha D^\rho_\beta
-\delta^\rho_\beta D^\sigma_\alpha )
$$
gives an invariant operator
$$
\Do^{\rho\sigma}_{\alpha\beta} : \ce_{\gamma\cdots \delta}^{\mu\cdots
  \nu}[w]\to \cs^{\rho\sigma}_{\alpha\beta}\otimes \ce_{\gamma\cdots \delta}^{\mu\cdots
  \nu}[w] .
$$
Similarly we define $\Do^{\rho\sigma\mu}_{\alpha\beta\gamma}$ by
$$
\Do^{\rho\sigma\mu}_{\alpha\beta\gamma}:=\frac{1}{3}( (D^\rho_\alpha
\Do^{\sigma\mu}_{\beta\gamma}+ D^\sigma_\beta
\Do^{\rho\mu}_{\alpha\gamma} + D^\mu_\gamma \Do^{\rho\sigma}_{\alpha\beta}
-\delta^\rho_\alpha\Do^{\sigma\mu}_{\beta\gamma}-\delta^\sigma_\beta
\Do^{\rho\mu}_{\alpha\gamma} 
-\delta^\mu_\gamma \Do^{\rho\sigma}_{\alpha\beta} )
$$
and so on for $\Do^{\alpha\cdots \delta}_{\rho\cdots
  \nu}$. 
Notice that the construction of these is designed in such a way that the
resulting operators are annihilated if composed (contracted) with
$Y_\nu^B$ on any index.

\smallskip
\noindent{\bf The Splitting Machinery.} 
In terms of the algebraic projectors and embeddings introduced in the last
section, the twistor-D  operator is given by
\begin{equation}\label{dform}
D^\rho_\alpha f= X^\rho_{R'}Y^A_{\alpha} \nd_A^{R'}f+ w
\xi^\rho_\alpha f ,
\end{equation}
where $f$ is any weighted twistor-spinor object.
Using this and the expressions (\ref{twconn}), (\ref{twconn2}) for the
twistor connection, the following identities are easily established:
\begin{equation}\label{diffids}
\begin{array}{lclclcl}
D^\rho_\alpha X^\beta_{C'}&=& X^\rho_{C'}\lambda^\beta_KY^K_\alpha&
      &D^\rho_\alpha Y^C_\beta&=& -Y^C_\alpha X^\rho_{K'}\xi^{K'}_\beta\\
D^\rho_\alpha \xi^{S'}_\beta &=& \Rho^{\rho S'}_{\alpha\beta} &
      &D^\rho_\alpha \lambda^\sigma_B& =&-\Rho^{\rho\sigma}_{\alpha B}\\
X^\alpha_{B'}D^\gamma_\alpha f &= &wX^\gamma_{B'} f & 
      &Y_\gamma^B D^\gamma_\alpha f &=& 0\\
\xi_\gamma^{B'} D^\gamma_\alpha f &=& D^{B'}_\alpha f&
&\lambda^\alpha_{B}D^{B'}_\alpha f &= &  \nd^{B'}_{B} f,
\end{array}
\end{equation}
where, again, $f$ is any weighted twistor-spinor and we write\\
$
\Rho^{\rho\sigma}_{\alpha\beta}:=
\Rho^{R'S'}_{AB}X_{R'}^\rho X_{S'}^\sigma Y^A_\alpha Y^B_\beta $,
$\Rho^{\rho S'}_{\alpha\beta}:=\Rho^{R'S'}_{AB}X_{R'}^\rho Y^A_\alpha
Y^B_\beta $, $\Rho^{\rho\sigma}_{\alpha B}:=
\Rho^{R'S'}_{AB}X_{R'}^\rho X_{S'}^\sigma Y^A_\alpha  $, etcetera.

Notice also that the objects $\xi_\alpha^{B'}$ \label{xistuff} and $
\lambda^\beta_A$ describing the splitting of the twistors can be
viewed as the projecting parts of
$\xi_\alpha^\beta:=\xi^{-1}D_\alpha^\beta \xi $ and
$\delta^\beta_\alpha-\xi^\beta_\alpha$, respectively.
   
\smallskip
\noindent {\bf $D$-Curvature.}
For $f\in\ce[w]$ the projecting part of $D_\alpha^\rho f$ is
$\sfrac{1}{p}X^\alpha_{P'} D_\alpha^{P'} f=wf$.  
Although this is 0th order in $f$, this part of $ D_\alpha^\rho f$
behaves like a first order operator because of the weight factor, $w$.
In particular $\sfrac{1}{p}X^\alpha_{P'} D_\alpha^{P'}$ satisfies a
Leibniz rule and so therefore so does $D^\rho_\alpha$. 
It follows immediately that, acting on $\ce^\mu [w]$,
$[D^\rho_\alpha ,D^\sigma_\beta]$ decomposes into a 0th order
curvature part and a 1st order torsion part. In fact it is an
elementary exercise using the identities \nn{skewrho1} and
\nn{diffids} 
to verify that 
\begin{equation}\label{Dcurv}
[D^\rho_\alpha ,D^\sigma_\beta] v^\mu  = 
W^{\rho\sigma\mu}_{\alpha\beta\gamma} v^\gamma 
-W^{\rho\sigma\nu}_{\alpha\beta\gamma} D_\nu^\gamma v^\mu
+\delta^\sigma_\alpha D^\rho_\beta v^\mu- 
\delta^\rho_\beta D^\sigma_\alpha v^\mu,
\end{equation}
where,
\begin{equation}\label{WW}
W^{\rho\sigma\mu}_{\alpha\beta\gamma} = X_{A'}^\rho X_{B'}^\sigma
Y^A_\alpha Y^B_\beta W^{A'B'\mu}_{AB\gamma} .
\end{equation}

\section{Invariant Theory}\label{invop}

Recall that each choice of  scale determines the linear connection 
$\nd^\xi$ on $\cg_0$. We shall write $\Gamma_{(\xi)}$ for the coefficients
of this connection $\nd^\xi$ in some coordinate frame. 
The linear connections
$\nd^{\xi}$ are clearly expressed through the normal Cartan connection
$\omega$ on $\cg$ and vice versa (this is one of the important aspects of the
Rho-tensor $\Rho_{ab}$, cf. Appendix \ref{A}). Thus we use the 
explicit definition of invariance given below. We use this approach for 
simplicity, but we should like to point out that there are more natural points
of view fitting nicely into the general concepts as developed in \cite{KMS}.
In particular, some of our polynomiality assumptions follow then
automatically.

\smallskip
\noindent {\bf Invariant operators.} Let $V$ and $ U$ be
finite dimensional $P$-modules with $ V$ irreducible. 
A {\em (coupled) invariant operator} 
on $\cv$ taking values in $\cu$ 
is a well defined differential operator $ \cv \to \cu$ which
may depend polynomially on the finite jets of the functions
$\Gamma_{(\xi)}$ as well
as polynomially on the finite jets of $\cv$ and which is independent
of the choice of local coordinate frame and scale $\xi$. By the very
definition, such an invariant must be intrinsic to the AG-structure and so, 
when evaluated, depends only on the section of
$\cv$ and the normal Cartan connection $\omega$ of the structure. 
It is clearly sufficient
to deal with the case that the invariant is homogeneous in $\cv$ and
we shall henceforth assume that coupled invariants are homogeneous in
this way. 

We say an invariant (and semi-invariants as described below)
has order $(\ell,m)$ if, in some scale $\xi$, it is:\\
(1) of order $\ell$ as
an operator on $ v\in \Gamma(\cv)$ and,\\
(2) in any coordinates, as an operator on
the functions $\Gamma_{(\xi)}$, it is of order $\geq m$ with equality
in some set of coordinates.\\ 
We will also
describe such an invariant as being of order $ k$ where $
k:=$max$(\ell,m)$. In the special case that the invariant is homogeneous of
degree 0 in the section $v$ then it is an invariant of the
structure. On the other hand if the invariant is of degree 1 in $ v$
then it is a linear invariant operator on $v$.

If the invariant takes values in $
\cu$ where $ U$ is an irreducible $P$-module
we will describe the invariant as an {\em irreducible invariant}. 
We may also restrict the definition of our operators to some subcategory of
the structures in question. For example we may require they are locally
flat, or torsion free, etc.

We will show below that the twistor-D operator is a universal
invariant differential operator in the sense that all 
coupled invariant operators arise in an appropriate sense from
concatenations of this operator and its curvature.

\smallskip\noindent{\bf Some Examples.}
Note that the invariance of the exterior derivative on functions is
implicit in the definition of the twistor-D operator. If $f\in \ce$
then $D^{P'}_\alpha f = (\nd^{P'}_A f\:\: 0)$ so the projecting part
of $D^\rho_\alpha f$ is $\nd^{P'}_A f$ and thus this is invariant. 
Similarly on 
$$
\left(\begin{array}{c} v^A \\ v^{A'}\end{array}\right)
=v^\alpha\in \ce^\alpha ,
$$
$D^{P'}_\alpha v^\beta =
(\nd^{P'}_A v^\beta \:\:
0)$  and so the projecting part of $D^\rho_\alpha
v^\beta$ is  $\nd^{P'}_A v^B+\delta^B_A v^{P'}$. The equation obtained
by setting this to zero is the usual twistor equation \cite{BaiE}. 

For a second order example consider $\Do^{\rho\sigma}_{\alpha\beta} f$
for $f\in \Gamma\ce[w]$ (or $f\in \Gamma\ce_{\alpha\cdots \gamma}^{\rho\cdots
  \mu}[w]$, with indices suppressed). It is easily established that
$$
\Do_{A\beta}^{P'\sigma} f=\left(\begin{array}{cc}
0&0\\
(\nd^{P'}_A\nd^{S'}_B) f + w{\sf S}^{P'S'}_{AB} f &
w\delta^{S'}_{B'}\nd^{P'}_A f -\delta^{P'}_{B'}\nd^{S'}_A f
\end{array}\right)
$$
and 
$$
\Do_{A'\beta}^{P'\sigma} f=\left(\begin{array}{cc}
0&0\\
w\delta^{P'}_{A'}\nd^{S'}_B f -\delta^{S'}_{A'}\nd^{P'}_B f &
w(w\delta^{P'}_{A'}\delta^{S'}_{B'}- \delta^{S'}_{A'}\delta^{P'}_{B'})f
\end{array}\right)
$$
and these two ``matrices'' display all the non-vanishing parts of
$\Do^{\rho\sigma}_{\alpha\beta} f$. 
Here, and below,  $(\nd_a\nd_b)$
means $\nd_{(a}\nd_{b)}$ and, recall, $ S_{ab}:=\Rho_{(ab)}$. If $w=0$
the projecting part of this is $\nd^{P'}_B f$. On the other hand if
$w=1$ then the projecting part of $\Do^{(\rho\sigma)}_{\alpha\beta} f$
is 
$$
(\nd^{P'}_{(A}\nd^{S'}_{B)}) f + {\sf S}^{P'S'}_{(AB)} f
$$
and so this is an invariant operator. Similarly if $w=-1$ then clearly
the projecting part of $\Do^{[\rho\sigma]}_{\alpha\beta} f$ is the
invariant operator
$$
(\nd^{P'}_{[A}\nd^{S'}_{B]}) f - {\sf S}^{P'S'}_{[AB]} f .
$$
In the case that $p=2$ we may contract this with $\vol_{P'S'}$ to
yield the invariant operator 
\begin{equation} \label{boxab}
\Box_{AB} f:= \vol_{P'S'}((\nd^{P'}_{A}\nd^{S'}_{B}) 
 - {\sf S}^{P'S'}_{AB} ) f .
\end{equation}
If also $q=2$ this is the usual conformally invariant Laplacian
or Yamabe operator $(\Delta-\sfrac{1}{6}R) f$ where $R$ is the Ricci
scalar curvature.

On the other hand if $w\neq -1,0,1$ then the projecting part of
both $\Do^{(\rho\sigma)}_{\alpha\beta} f$ and
$\Do^{[\rho\sigma]}_{\alpha\beta} f$ is a non-zero multiple of $f$.
In fact it is an easy consequence of this observation and the theorem
\ref{main} that, on weighted
functions of weight $w\neq -1,0,1$, there are no linear
invariant operators of order $\leq 2$ which are non-trivial on flat 
structures.

Note that in, for example, the $w=-1$ case above the operator 
$(\nd^{P'}_{[A}\nd^{S'}_{B]}) f - {\sf S}^{P'S'}_{[AB]} f$ may be
described explicitly by the formula
$$
\xi_\rho^{P'}\xi_\sigma^{S'}\lambda_A^\alpha\lambda_B^\beta 
\Do^{[\rho\sigma]}_{\alpha\beta} f .
$$
It is useful to think of this as a composition of 
$$
\Do^{[\rho\sigma]}_{\alpha\beta}:\ce[-1]  \to {\cal Q}^{\rho\sigma}_{\alpha\beta}  
$$
with 
\begin{equation}\label{comp1}
\xi_\rho^{P'}\xi_\sigma^{S'}\lambda_A^\alpha\lambda^\beta_{B}:
{\cal Q}^{\rho\sigma}_{\alpha\beta} \to
\ce^{[P'S']}_{[AB]}[-1] .
\end{equation}
Here ${\cal Q}^{\rho\sigma}_{\alpha\beta} $ is the minimal natural 
sub-bundle of $ \cs^{\rho\sigma}_{\alpha\beta}[-1]$ which contains the
image of $ \Do^{[\rho\sigma]}_{\alpha\beta}$ on $ \ce[-1]$. This is
induced by a $ P$-submodule, say $ H$, of the representation inducing  $
\cs^{\rho\sigma}_{\alpha\beta}[-1]$ and the invariant map \nn{comp1} arises
from a $P$-homomorphism from $H$ to $U$, where $ U$ is the $ P$-module
inducing $\ce^{[P'S']}_{[AB]}[-1] $. It is clear that one can use such
$P$-homomorphisms composed with the $ \Dok$ operators to proliferate
invariants. The content of theorem \ref{main} is that all invariants
arise  this way.

\smallskip\noindent{\bf An easy proposition.}
We will observe here
that via the twistor-D operator we obtain
a special description of the jet bundle associated to twistor sub-bundles.

Let $\cv$ be any subbundle of a weighted twistor bundle and let us
write $\Dok$ for the linear differential operator
$$
\Dok : \cv \to (\ce\oplus \cs^\rho_\alpha \oplus
\cs^{\rho\sigma}_{\alpha\beta} \oplus \cdots \oplus
\cs^{\mbox{\tiny $\rho\cdots\mu$}}\sub{\alpha\cdots \gamma}{k})\otimes \cv
$$
given by
$$
f\mapsto f\oplus \Do^\rho_\alpha f \oplus
\Do^{\rho\sigma}_{\alpha\beta}f \oplus \cdots
\oplus \Do^{\rho\cdots\mu}_{\alpha\cdots \gamma} f .
$$
  Recall that any $k$th order differential operator on a bundle factors
  through the associated bundle of $k$-jets.  That is, any $k$th order
  linear invariant differential operator, taking values in a bundle
  $\cu$, $\cv \to \cu$, is equivalent to a bundle morphism
  $J^k(\cv)\to \cu$. In particular, $\Do^{(k)}$ factors through a linear
mapping on the $k$th jet prolongation $J^k(\cv)$. The image of $\Do^{(k)}$
fills a vector sub-bundle of 
$(\ce\oplus \cs^\rho_\alpha \oplus
\cs^{\rho\sigma}_{\alpha\beta} \oplus \cdots \oplus
\cs^{\rho\cdots\mu}_{\alpha\cdots \gamma})\otimes \cv$, which we denote by 
$\cjk(\cv)$. 

\begin{proposition} \label{bunjets}
Let $\cv$ be any subbundle of a weighted twistor bundle of weight $w$.
The operator $ \Dok$ determines a bundle isomorphism, 
$$
J^k(\cv)\cong\cjk(\cv).
$$
\end{proposition}
\begin{proof}
 In view of the definition of $\cjk(\cv)$, the
  operator $\Dok$ clearly determines a bundle epimorphism $
  J^k(\cv)\to\cjk(\cv)$. 
That this is also injective follows by counting
  dimensions: Consider $f\in\ce[w]$. Observe that the injecting part
  of
$$
\Do^{\mbox{\tiny $\rho\cdots\mu$}}\sub{\alpha\cdots \gamma}{k} f
$$
is of the form
\begin{equation}\label{form}
\underbrace{\nd_{(a}\cdots \nd_{d)}}_{k} f + ~(\mbox{lower order terms}) .
\end{equation}
All other parts of $ \Do^{\rho\cdots\mu}_{\alpha\cdots \gamma} f $ are
of order at most $k-1$ and so, by repeated use of \nn{form}, can be
expressed polynomially in terms of
$$
\Do^{\mbox{\tiny $\rho\cdots\sigma$}}\sub{\alpha\cdots
  \beta}{\ell} f
$$
for $\ell\leq k-1$. Thus 
$$
\cjk(\ce[w])/\cj^{(k-1)}(\ce[w])\cong (\odot^k \ce_a)\otimes \ce[w]
$$
but $(\odot^k \ce_a)\otimes \ce[w]\cong
J^k(\ce[w])/J^{k-1}(\ce[w])$. In fact it is easily seen that, by an
almost identical argument, we have the more general result,
$$
\cjk(\cv)/\cj^{(k-1)}(\cv)\cong (\odot^k \ce_a)\otimes \cv 
\cong J^k(\cv)/J^{k-1}(\cv)
$$
and so, by induction on $k$, the fibre dimension
of $\cjk(\cv)$ is the same as the fibre dimension of $J^k(\cv)$.
\end{proof}

Although the proposition above is inspiring we need to consider
slightly more general structures to obtain all invariants. These are
defined above theorem \ref{main}.

In the meantime we need to understand the 
general invariants of $\nd_{\xi}$, viewed as 
affine connections. 
In order to distinguish them from the 
AG-invariants, we will call these {\em semi-invariants} in the sequel.

\smallskip\noindent{\bf Semi-invariants and their normal form.}
As above, let $V$ and $ U$ be
finite dimensional $P$-modules, with $ V$-irreducible, and let $\cv$ and
 $\cu$ be the corresponding natural bundles. A {\em coupled
 semi-invariant operator} (which we will often abbreviate to 
{\em semi-invariant}) on $\cv$ taking values in $\cu$ is a universal
formula which is polynomial in the coordinate derivatives of the
functions $\Gamma_{(\xi)}$ and coordinate derivatives of the
components of $\cv$ (in some local frame) which is
independent of the choice of local coordinates and frame (but may not be
 independent of the choice of scale $\xi$).
Thus, for each choice of scale $\xi$, a semi-invariant is a
differential operator $ \cv \to \cu$ which may depend
polynomially on the finite jets of the functions $\Gamma_{(\xi)}$ as
well as polynomially on the finite jets of $\cv$.
 Note that a coupled
 invariant operator is a coupled
 semi-invariant operator which, in addition, is independent of the
 choice of scale $\xi$. As for invariants, semi-invariants will be
 deemed {\em irreducible} if they take values in irreducible natural bundle. 
 
 It is easy to write down some examples of such semi-invariants. The
 curvature $R_{ab}^\vb{}^c_d$ and torsion $T_{ab}{}^c$ of $\nd^\xi$ are
 polynomial in the finite jets of the functions $\Gamma_{(\xi)}$ and
 it is a classical result that these objects are tensorial and so
 are semi-invariants. Thus the irreducible parts
 of these tensors are irreducible semi-invariants. The 
 objects $F_{ab}{}^c$ and $U_{ab}^\vb{}^C_D$ are examples. In fact due to
 the invariance of the covariant derivative it is easily verified that
 any contraction involving
 covariant derivatives of $v\in\Gamma(\cv)$ and covariant derivatives
 of the curvature $R_{ab}^\vb{}^c_d$ and the torsion $T_{ab}^c$ is 
 a semi-invariant. 
 For example
$$
v^a(\nd_a v^c)(\nd_cU_{de}^\vb{}^H_I)U_{fg}^\vb{}_H^I 
$$
is a semi-invariant.
 We will write
 $\contr(\nd^\xi,T,R,v)$ to symbolically indicate such  contractions.
 We will observe that, in fact, all semi-invariants arise this way and
 this leads to a standard way of expressing semi-invariants.
 
 Let us fix a scale $\xi$. Note that it follows easily from
 proposition \ref{Taylor} and proposition \ref{normframe} 
that a semi-invariant may be expressed as
 a polynomial in the components of the covariant derivatives of the
 torsion and curvature of $\nd^{(\xi)}$ and the  components of the covariant
 derivatives of the section $v\in\Gamma(\cv)$.  At each point of the
 manifold a   semi-invariant is a polynomial in the components of
 these tensors (that is the list of tensors which give, at that point,
 the various covariant derivatives of $T_{ab}{}^c$, $R_{ab}^\vb{}^c_d$ and
 $v$
 to the required order) which is covariant under the action of $
 \SL(p)\times{\SL}(q)$. Thus it follows from the complete reducibility
 of finite dimensional $(\SL(p)\times\SL(q))$-modules and Weyl's
classical invariant theory \cite{weyl} that {\em any} such
 semi-invariant can be expressed as a linear combination of basic
 semi-invariants of the form $\contr(\nd^\xi,T,R,v)$ as claimed.

 Consider then a semi-invariant expressed as a linear combination of
 contractions $\contr(\nd^\xi,T,R,v)$. First observe that by
 substituting for $R_{ab}^\vb{}^c_d$ using formulae \nn{curvdec} 
and \nn{Wsplit}
 we see that our typical semi-invariant may be
 re-expressed in terms of covariant derivatives of the objects $v$
 (with indices suppressed), $T_{ab}{}^c$, $U_{ab}^\vb{}^C_D$,
 $U_{ab}^\vb{}^{C'}_{D'}$ and $\Rho_{ab}$. We might write
 $\contr(\nd^\xi,T,U,U',\Rho,v)$ for the basic terms of this new
 expression. Finally we observe that the semi-invariant can be written
 as described in the following lemma which we regard as a {\em normal
 form} for semi-invariants.
\begin{lemma}\label{seminormf}
  Any semi-invariant of order $k$ may be expressed as a linear
  combination of contractions involving the tensors
$S_{a\cdots d}\in \odot^{m}\ce_a$  
for $0 \leq m \leq k$, and various covariant derivatives of the
  objects $T_{ab}{}^c$, $U_{ab}^\vb{}^C_D$, $U_{ab}^\vb{}^{C'}_{D'}$,
  $Q_{abc}$ and $v\in\Gamma(\cv)$.
\end{lemma}
\begin{proof}
Recall that 
$$
Q_{abc}= \nd_{[a}\Rho_{b]c}-\frac{1}{2}T_{ab}{}^e\Rho_{ec},
$$
(as in \nn{twcurv}).
Note that it is easily verified, by considering possible Young
symmetrizers and using \nn{skewrho}, that 
the $(m-2)$nd covariant derivative of $\Rho_{ab}$,
$$
\underbrace{\nd_{a}\cdots \nd_b}_{m-2}\Rho_{cd},
$$
may, {\em up to lower order terms} which involve covariant derivatives of
the curvature and torsion, be expressed as a linear combination  of
the tensors (cf. (\ref{skewrho}))
$$
S\sub{a\cdots d}{m},\quad
\underbrace{\nd_{a}\cdots \nd_b}_{m-3}\nd_{[c}\Rho_{d]e} ~~\mbox{ and
  }
\underbrace{\nd_{a}\cdots \nd_b\nd_e}_{m-1}T_{cd}{}^e
$$
 Thus, by replacing $\nd_{[c}\Rho_{d]e}$
with $ Q_{cde}+ \sfrac{1}{2}T_{cd}{}^f\Rho_{fe}$ it is clear that
the tensors $\nd_{a}\cdots \nd_b\Rho_{cd}$ may, up to lower order terms
which involve covariant derivatives of the curvature and torsion, be
expressed as a linear combination of the tensors 
$$
S\sub{a\cdots d}{m},\quad
\underbrace{\nd_{a}\cdots \nd_b}_{m-3}Q_{cde} ~~\mbox{ and
  }
\underbrace{\nd_{a}\cdots \nd_b\nd_e}_{m-1}T_{cd}{}^e .
$$
The lemma follows by first expressing the semi-invariant in the manner
last described above and then
repeatedly using this observation to replace all occurrences of
covariant derivatives of $\Rho_{ab}$, starting with the highest
order. 
\end{proof}

\smallskip\noindent{\bf The Main theorem.}
Note that the covariant derivatives of tensors and spinors can be
expressed in terms of components of the covariant derivatives of
$\Rho_{ab}$ and components of  the twistor-D operator acting on 
appropriate twistors
via the machinery of section \ref{twistcalc}. For example consider
$\nd_a v_B$ where $v_B\in \ce_B[w]$. Let $ v_\beta:=Y_\beta^B
v_B$ and we have 
$$
\nd_a v_B= \xi^{A'}_\rho \lambda^\alpha_A D^\rho_\alpha
\lambda_B^\beta v_\beta .
$$
We now bring the $\lambda_B^\beta$ to the left of the twistor-D
operator using the appropriate identities from \nn{diffids} and \nn{algids}. 
We obtain
$$
\nd_a v_B= \xi^{A'}_\rho \lambda^\alpha_A\lambda_B^\beta (D^\rho_\alpha
v_\beta).
$$
Thus
\begin{align*} 
\nd_a\nd_b v_C =\ & \xi^{A'}_\rho \lambda^\alpha_A D^\rho_\alpha
\xi^{B'}_\sigma \lambda^\beta_B D^\sigma_\beta
\lambda^\gamma_C v_\gamma\\
=\ &\xi^{A'}_\rho \lambda^\alpha_A D^\rho_\alpha (\xi^{B'}_\sigma
\lambda^\beta_B \lambda^\gamma_C (D^\sigma_\beta v_\gamma))\\
=\ &\xi^{A'}_\rho \lambda^\alpha_A \xi^{B'}_\sigma
\lambda^\beta_B \lambda^\gamma_C (D^\rho_\alpha D^\sigma_\beta
v_\gamma)
-\Rho^{A'C'}_{AB} X^\beta_{C'} \xi^{B'}_\sigma \lambda^\gamma_C
(D^\sigma_\beta v_\gamma)
\\&-\Rho^{A'C'}_{AC} X^\gamma_{C'}\xi^{B'}_\sigma \lambda^\beta_B
(D^\sigma_\beta v_\gamma) 
\end{align*}
Continuing
in this fashion it is easily seen that 
$$
\underbrace{\nd_a\cdots \nd_b}_\ell v_B
$$
may be expressed in terms of components of  
$$
\underbrace{D^\rho_\alpha\cdots D^\sigma_\beta}_\ell v_\gamma
$$
and lower order terms which polynomially involve the components of  
$$
\underbrace{D^\rho_\alpha\cdots D^\sigma_\beta}_m v_\gamma ,
$$
for $m\leq \ell-1$, and the components of covariant derivatives, to
order $\ell-2$, of $\Rho_{ab}$.

These observations lead us to the next lemma which is the key to the 
proof of the theorem in this section. As before let $\cv$ be an
irreducible natural bundle and recall that we may identify this  with
a twistor sub-bundle
(lemma \ref{spintw}). Since we are suppressing the
indices on the section $v\in\cv$ we will write
$\tilde{v}\in\tilde{\cv}$ 
for the
corresponding section of the appropriate twistor bundle. The section
$v$ is recovered explicitly by contracting $\tilde{v}$ with the
projectors $\xi^{A'}_\alpha$ and $\lambda^\alpha_A$. (For example, 
in the example just above $v=v_B\in\ce_B[w]$ and
$\tilde{v}=v_\beta\in\ce_\beta[w]$ with $ v_B=\lambda_B^\beta v_\beta$.) 
\begin{lemma}\label{normf}
A coupled invariant differential operator $I$ of order $(\ell,m)$ may
be expressed as a universal polynomial expression in the
components of $\Do^{(\ell)} \tilde{v}$ and $\Do^{(k')} W$ where
$k'={\rm max}(\ell-1,m)$ .
\end{lemma}
In this lemma, and henceforth,  $\Do^{(m)} W$ means $\Do^{(m)}$
applied to $W^{\rho\sigma\gamma}_{\alpha\beta\delta}$. This is to be
distinguished from $(D)^{(m)} W$ which we will use to mean simply an $
m$-fold application of the twistor-D operator to
$W^{\rho\sigma\gamma}_{\alpha\beta\delta}$. 
\begin{proof}
  We may suppose that at first we have chosen a scale $\xi$ and the
  invariant is expressed in normal form as in lemma \ref{seminormf}.
  We will first observe that, in rewriting this expression, covariant
  derivatives of $T_{ab}{}^c$, $U_{ab}^\vb{}^C_D$, $U_{ab}^\vb{}^{C'}_{D'}$,
  and $Q_{abc}$ may be eliminated in favour of components of
  $\Do^{(m)} W$ and lower order terms and similarly covariant
  derivatives of $v\in\Gamma(\cv)$ may be 
  
  Note that each of the objects
  $T_{ab}{}^c$, $U_{ab}^\vb{}^C_D$, $U_{ab}^\vb{}^{C'}_{D'}$, and $Q_{abc}$
  may be obtained linearly from
  $W^{\rho\sigma\gamma}_{\alpha\beta\delta}$ via the projectors
$X^\alpha_{A'}, Y_\alpha^A, \xi_\alpha^{A'}, \lambda^\alpha_A$. For
  example
$$
U_{ab}^\vb{}^C_D= \xi^{A'}_\rho \lambda_A^\alpha \xi^{B'}_\sigma \lambda_B^\beta
Y^C_\gamma \lambda_D^\delta
W^{\rho\sigma\gamma}_{\alpha\beta\delta} .
$$
We re-express the invariant as follows. We make the substitutions
for each of $T_{ab}{}^c$, $U_{ab}^\vb{}^C_D$, $U_{ab}^\vb{}^{C'}_{D'},
Q_{abc}$ and $v$ in terms of components of
$W^{\rho\sigma\gamma}_{\alpha\beta\delta}$ and $\tilde{v}$, and we
replace each $\nabla_a$ with $\xi^{A'}_\rho \lambda^\alpha_A
D^\rho_\alpha$. Next we further re-express by moving each
$X^\alpha_{A'}, Y_\alpha^A, \xi_\alpha^{A'}$ and $\lambda^\alpha_A$ to
the left of any twistor-D operators. The new expression for the
invariant involves components of concatenations of twistor-D operators
acting on $W^{\rho\sigma\gamma}_{\alpha\beta\delta}$ and $\tilde{v}$
and covariant derivatives of $\Rho_{ab}$ and various valence
$S$-tensors.  These covariant derivatives of $\Rho_{ab}$ all turn up
via the identities \nn{diffids}. From this observation it is
immediately clear that the order of any of these covariant derivatives
of $\Rho_{ab}$ is strictly less than $k:=$max$(\ell,m)$.  In fact, by
elementary representation theory arguments, one can show that the
order of any of these covariant derivatives of $\Rho_{ab}$ is $\leq
\ell-2$ if $\ell>m$, and is $ \leq m-2$ otherwise. 
Now we replace, in
the last expression for the invariant, each maximal order $\nd_a\cdots
\nd_c\Rho_{de}$ with its expression in terms of the tensors $ S_{a
  \cdots e}$, $ \nd_a\cdots \nd_b Q_{cde}$, $\nd_a\cdots \nd_c \nd_f
T_{de}{}^f$, their transposes
and lower order terms.  
Next we replace each occurrence of
$\nd_a\cdots \nd_b Q_{cde}$ and $\nd_a\cdots \nd_c \nd_f T_{de}{}^f$
with their expressions in terms of components of
$(D)^{k'} W$ and lower order covariant
derivatives of $\Rho_{ab}$. Continuing in this fashion it is clear
that finally we are left with an expression involving only $S$-tensors
and the components of concatenations of the twistor-D operator on
$W^{\rho\sigma\gamma}_{\alpha\beta\delta}$ and $\tilde{v}$.  It is
easily seen using \nn{Dcurv} that this may be re-expressed in terms of
components of $\Do^{(\ell)} \tilde{v}$, $ \Do^{(m)} W$ and the
components of the $S$-tensors.

Now let us write the invariant $I$ as a sum of two parts
$$
I=A+B
$$
where the part $A$ consists of all terms which involve no
components of the $ S$-tensors while $B$ is the remaining part which
consists of all terms which do involve the $ S$-tensors. Let us
choose a point $q$ and consider changing the scale of $\xi$ by a
factor $\Omega$ so that $\Upsilon_a(q)=0$. Under such a transformation
it is clear that, at $ q$, the $ A$ part of $ I$ is invariant as the
transformation of the components of an invariant twistor depends only
on the first derivative of $ \Omega$. On the other hand $ I$ is
invariant under any transformation. Thus it follows that, under 
transformations such that $\Upsilon_a(q)=0$, $ B$ must also be
invariant. But on the other hand $B$ vanishes in a normal scale $\xi_q$
(see \nn{normsc} in section \ref{normal}) since in this scale all
the $S$-tensors vanish at $ q$. As observed in remark \ref{noneedupq},
such a scale can be achieved by a
transformation with $\up_a(q)=0$. Thus $ B$ must vanish at $ q$. Since
we may perform this calculation at any point it follows that $ B$ 
vanishes everywhere so  $I=A$ and the lemma is proved. 
\end{proof}

Before we can discuss the main theorem we will need some special
notation. Recall that $\cj^\ell(\cv)$ was defined to be the subbundle of the
natural bundle 
$$
\cv^{(\ell)}:=(\ce\oplus \cs^\rho_\alpha \oplus
\cs^{\rho\sigma}_{\alpha\beta} \oplus \cdots \oplus
\cs^{\mbox{\tiny $\rho\cdots\mu$}}\sub{\alpha\cdots \gamma}{\ell})\otimes
\tilde{\cv}
$$
determined by the image of the 
invariant operator $\Do^{(\ell)}$ on $\tilde{\cv}$. Note that in
general $\cj^\ell(\cv)$ will not itself be  a natural bundle as the
algebraic properties of its fibres vary over $M$. Suppose that $\mu$ and
$V^{(\ell)}$ denote, respectively, the $P$-representation and 
representation space 
inducing the natural bundle $ \cv^{(\ell)}$ displayed.
If $\tilde{v}$ is a section of $\tilde{\cv}$ then $\Do^{(\ell)}(\tilde{v})$ 
is a section of
$\cj^\ell(\cv)$. That is $\Do^{(\ell)}(\tilde{v})$ is a function
$$
\Do^{(\ell)}(\tilde{v}) : \cg \to V^{(\ell)}
$$ 
which is homogeneous, $\Do^{(\ell)}(v)(x.p)=\mu(p^{-1})\Do^{(\ell)}(v)(x)$, $x\in
\cg$ and $p\in P$. Note that in general this function is not
surjective.  Note also that the image of $\Do^{(\ell)}$ depends on the
underlying structure of $M$, that is on the normal Cartan bundle equipped by
the normal Cartan connection. 

Our next step is to construct a sort of smallest natural bundle
which could accommodate the values of $\Do^{(\ell)}$.
Let us fix a point $q\in M$ and a coordinate neighbourhood $Q={\Bbb R}^{pq}$
centred at $q$ (in fact we may forget about
our manifold $M$ and we work just over ${\Bbb R}^{pq}$ for the while).
Let us consider all possible normal Cartan connections on the trivial Cartan
bundle $\cg=Q\times P$ and for each such normal Cartan connection $\omega$, 
let us write 
$$
\cj^\ell_o(\cv,\omega)
$$
to denote the span of the image of $\Do^{(\ell)}(\tilde{v})$, on the fibre of 
$\cg$ over $q$, as we vary over all possible argument sections $v$. Note
that $\cj^\ell_o(\cv,\omega) $ is a well defined $P$-submodule of
$V^{(\ell)}$. Now let 
$$
\cj^\ell_o(\cv):=\langle\cup_\omega \cj^\ell_o(\cv,\omega)\rangle
,$$
the span of the union which is taken 
over all possible normal Cartan connections on $\cg\to Q$.
Then $\cj^\ell_o(\cv)$ is also a well defined $P$-submodule
of $V^{(\ell)}$ and the corresponding natural subbundle in $\cv^{(\ell)}$ 
is the smallest one 
containing all possible subbundles $\cj^\ell(\cv)$. 

Next we observe that we can consider a similar `generic natural bundle' 
for the curvature. Write $W^{(m)}$ for the $P$-module inducing
$$
(\ce\oplus \cs^\rho_\alpha \oplus
\cs^{\rho\sigma}_{\alpha\beta} \oplus \cdots \oplus
\cs^{\mbox{\tiny $\rho\cdots\mu$}}\sub{\alpha\cdots
  \gamma}{m})\otimes \ce^{\upsilon\phi\zeta}_{\delta\epsilon\vartheta} .
$$
Then, for each normal Cartan connection,  $\Do^{(m)} W$ takes values in
$W^{(m)}$ and the span of the image of 
$\Do^{(m)} W$ on the
fibre of $\cg$ over $ q\in M$ is a $ P$-submodule
of  $W^{(m)}$ that we will denote $\cj^m_o(\cw,\omega)$. 
In analogy with the above we let
$$ 
\cj^m_o(\cw):=\langle\cup_\omega \cj^m_o(\cw,\omega) \rangle
$$ 
here the union is taken over all possible normal Cartan connections
defined locally at the fixed point $q\in M$.
This is clearly a $P$-submodule of $W^{(m)}$.
Now write
$$
\cj^{m,t}_o(\cw):=\oplus_{i=0}^t  \odot^i \cj_o^m(\cw).
$$

What we really need is a `generic fibre' for a mixed case. Write 
$$
\cj^{\ell,s,m,t}_o(\cv,\cw)
$$ 
for the $P$-submodule of 
$(\odot^s V^{(\ell)})\otimes (\oplus_{i=0}^t \odot^i W^{(m)})$  determined
by the image of 
$(\odot^s (\Do^{(\ell)} v)) \otimes 
(\oplus_{i=0}^t \odot^i \Do^{(m)} W)$, at $ q
\in M$, as we vary the normal Cartan
connection and for each such connection vary $\tilde{v}$ over all possible
sections of $\tilde{\cv}$. Note that $ \cj^{\ell,s,m,t}_o(\cv,\cw)$ is clearly
a $P$-submodule of $(\odot^s\cj^\ell_o(\cv))\otimes \cj^{m,t}_o(\cw)$.
Note also that the previous `generic fibres' are special cases of
this,  
$\cj^{\ell,s,m,0}_o(\cv,\cw)= \odot^s\cj^\ell_o(\cv)$ and
$\cj^{\ell,0,m,t}_o(\cv,\cw)=\cj^{m,t}_o(\cw)$. 

The following is the main theorem of this section. 
\begin{theorem}\label{main}  
  A coupled invariant operator $I$ which is homogeneous of degree $s$ 
on an irreducible bundle $\cv$, and 
  taking values in the natural bundle $\cu$, 
is equivalent to a $P$-homomorphism
$$
I_o: \cj^{\ell,s,m,t}_o(\cv,\cw) \to U . 
$$
That is there is a $1-1$ correspondence between such invariant
operators $ I$ and homomorphisms $ I_o$ as indicated.
\end{theorem}
\begin{proof}
$\Rightarrow$: Let $I(v,\omega)$ indicate the invariant $I$ evaluated
on a section $v$ of $ \cv$ and some particular normal Cartan
connection. Then $I(v,\omega): \cg\to U$ with the homogeneity property
$I(v,\omega)(x.p)=\sigma(p^{-1})I(v,\omega)(x)$ where $\sigma$ denotes
the inducing representation of $P$ on $U$. 

Lemma \ref{normf}, combined with standard polarisation techniques,
implies that there is a linear mapping $\tilde I$,
defined on the whole $P$-module 
$\odot^s V^{(\ell)}\otimes(\oplus_{i=0}^t \odot^i W^{(m)})$ 
(notice that our operator is
homogeneous in the arguments from $\cv$) such that
\begin{equation}\label{equiv}
I(v,\omega)(x) = \tilde I(
(\odot^s (\Do^{(\ell)} \tilde{v})) \otimes (\oplus_{i=0}^t \odot^i \Do^{(m)}
W(\omega))(x))
\end{equation}
for all $x\in\cg$. The mapping 
$\Phi(v,\omega):\cg\to \cj^{\ell,s,m,t}_o(\cv)$,
$$
x\mapsto (\odot^s (\Do^{(\ell)} \tilde{v})) \otimes (\oplus_{i=0}^t \odot^i
\Do^{(m)} W(\omega))(x)
$$ 
is $P$-equivariant too, and a general element in
$\cj^{\ell,s,m,t}_o(\cv)$ is a finite linear combination
$\sum_jc_j\Phi(v_j,\omega_j)$. Since $\tilde I$ is linear, the equivariance of
the compositions $I(v_j,\omega_j)=\tilde I\circ\Phi(v_j,\omega_j)$ implies that
the restriction $I_o$ of $\tilde I$ to $\cj^{\ell,s,m,t}_o(\cv)$ 
is $P$-equivariant, as required. This shows that all invariants $I$
arise from a $P$-homomorphism as in the theorem. 

\medskip

$\Leftarrow:$
The composition of $I_o$ with 
$(v,\omega)\mapsto (\odot^s (\Do^{(\ell)} \tilde{v})) \otimes 
(\oplus_{i=0}^t \odot^i \Do^{(m)} W)$, for each $v\in \Gamma(\cv)$ and 
normal Cartan connection $\omega$, 
is clearly 
a coupled invariant operator.
This shows that all invariants $I$
arise from a $P$-homomorphism as in the theorem. 
We complete the proof by showing that if $I_o\neq 0$ then $I\neq 0$.   
If $I_o\neq 0$ then  there exists an element of
$\cj^{\ell,s,m,t}_o(\cv,\cw)$ such that $I_o$ does not kill this element.
As mentioned above a general element of
$\cj^{\ell,s,m,t}_o(\cv,\cw)$ may be expressed by a finite linear
combination
$\sum_jc_j\Phi(v_j,\omega_j)$. Let this finite linear combination
represent, in particular, the element not killed by $I_o$. It follows that
one of the $\Phi(v_j,\omega_j)$ is not killed by $I_o$. That is there
exists a section $v=v_j$ and a Cartan connection $\omega=\omega_j$ such that 
$I_o(\Phi(v,\omega))\neq 0$. But this means the invariant differential
operator which is $I_o$ composed with $\Phi$
is non-trivial. Thus we have shown that any non-trivial $P$-hom yields a
non-trivial invariant operator as required.
\end{proof}


\section{New Invariant Operators}\label{newops}

As we have discussed, the case $p=q=2$ corresponds to the usual
four-dimensional conformal structures. 
Here we restrict attention to torsion-free AG-structures with $p=2$, $q>2$.
The main result is theorem \ref{weaks} which, for these geometries, 
gives curved analogues
for all the non-standard operators between differential forms.
We will deal with $q>2$ odd as well as even, but we would like to
point out that the cases of $q$ even are of particular interest as these
include all the quaternionic geometries. That is, when $q$ is even our
formulae and results below describe invariant operators on the 
quaternionic geometries. This is in the spirit of our simultaneous
treatment of AG-structures, so of course the formulae also give invariant
operators for the other geometries (i.e. those corresponding to the `real
split form' SL$(p+q,{\Bbb R})$) and we obtain similar
operators when $q$ is odd.  

At this point it is worthwhile to review the examples exposed above and
in particular the operator $\Box_{AB} f$ as displayed in \nn{boxab}.  Although this second order 
invariant operator does
not operate between forms it is closely related to the fourth order operators
we construct below.

To describe the operators it is useful to have some efficient and
concise notation for the bundles concerned. For this we will use Young
diagrams \cite{FH,ot}.
We will use these to indicate projections onto
irreducible representations of SL$(m)$. In our case we will in
particular use these for representations of SL$(p+q)$, representations
of ${\rm SL}(p)\times {\rm SL}(q)$, which are trivial with respect to
the SL$(p) $ factor, and the bundles these induce.
(Here, as usual, SL$(r)$ can mean either  SL$(r,{\Bbb R})$ or 
SL$(r,{\Bbb C})$ depending on which structures we are considering. The
comments here apply equally to both cases.) For example we
could write $\yiiI(\otimes^2\ce_A)$ or $ \yiiI\ce_{AB}$ to mean
$\ce_{(AB)}$. In fact we will shorten this notation further and simply
write $\yiiI \ce_A$ for this, that is $\yiiI\ce_A=\ce_{(AB)}$.  In this
notation the total number of boxes in the given Young diagram
indicates the required tensor power of the bundle. For diagrams of height and
width greater than 1 we adopt the convention that we symmetrized over
sets of indices corresponding to rows of the diagram first and then
with the result skew over sets of indices corresponding to the columns
of the diagram. For instance suppose we start with some general valence 3
spinor $A_{EFG}\in \ce_{EFG}[w]$.  If we first symmetrized over the
last two indices to form $B_{EFG}:=A_{E(FG)}$ and then on this result
skew on the first two indices to obtain $C_{EFG}:= B_{[EF]G}$, then
$$
C_{EFG}\in (\yiIIi\ce_F)[w].
$$
Although, for the sake of being concrete, we will suppose that this is
the convention adopted, nothing we do actually depends on this choice
of convention.

We will use this notation immediately in the construction of a special
invariant operator.  Recall that if $\ce^\ast[w]$ is any weighted
twistor bundle then we have the invariant operator
$\Do_{\alpha\beta\gamma\delta}^{\rho\sigma\mu\nu}:\ce^\ast[w]\to
\cf_{\alpha\beta\gamma\delta}^{\rho\sigma\mu\nu}\otimes\ce^\ast[w]$.
Equivalently we may view this as an operator
$$
D_{\alpha\beta\gamma\delta}^{R'S'U'V'}: \ce^\ast[w]\to
\ce_{\alpha\beta\gamma\delta}^{R'S'U'V'}\otimes \ce^\ast[w] .
$$
Thus there is an invariant operator
$$
D_{\alpha\beta\gamma\delta}:=\ep_{R'S'}\ep_{U'V'} D_{\alpha\beta\gamma\delta}^{R'S'U'V'}: \ce^\ast[w]\to
\ce_{\alpha\beta\gamma\delta} \otimes \ce^\ast[w-2] .
$$
Note that $ D_{\alpha\beta\gamma\delta}$ inherits some symmetry from
$\Do_{\alpha\beta\gamma\delta}^{\rho\sigma\mu\nu}$, in particular
observe that $ D_{\alpha\beta\gamma\delta}= D_{[\alpha\beta][\gamma\delta]}$.
For any $0\leq k\leq q-2$ and weight $ w\in {\Bbb R}$, let us 
write $\Box_{\alpha\beta\gamma\delta} $ for the non-trivial composition
of 
$$
{D_{\alpha\beta\gamma\delta}}:
(\yii{k}\ce_\sigma)[w] \to (\ce_{\alpha\beta\gamma\delta}\otimes
(\yii{k}\ce_\sigma)[w-2])
$$
with a Young projection
$$
(\ce_{\alpha\beta\gamma\delta} \otimes (\yii{k}\ce_\sigma)[w-2]\to 
(\wyii{k+2}\ce_\alpha)[w-2].
$$
This is clearly invariant for all $w$. Note also that it is an
elementary exercise to verify that there is such a composition which is
non-trivial and that it is unique up to a natural isomorphism of the
image bundle. (For example, in the $k=0$ case the main point is to
observe that $\yiiII$ turns up precisely once in the product $
\yiII\otimes\yiII$.) 


Before we state the theorem let us introduce one further item of
notation.  Let us write $\ch_\alpha$ for the subbundle of $\ce_\alpha$
which is naturally isomorphic to $\ce_A$ (c.f.\ $ \cf^\alpha$ of
section \ref{twistcalc}).  Here is the main result of this section.
\begin{theorem}\label{weaks}
Let $M$ be a torsion-free AG-structure, $p=2$, $q>2$. 
For each integer $k$ such that $0\leq k\leq q-2$ there is a fourth order 
invariant operator,
$$
\Box_{ABCD}:(\yii{k}\ce_E)[-k]\to (\wyii{k+2}\ce_E)[-k-2], 
$$
which is non-trivial on flat structures.

For each $k$ the operator is given by 
$$
{\Box_{\alpha\beta\gamma\delta}}:
(\yii{k}\ch_\vep)[-k] \to (\wyii{k+2}\ch_\vep)[-k-2].
$$
\end{theorem}

Before entering the proof of this theorem, we shall discuss 
the corresponding operators on the locally flat AG-structures since
their existence is a key to our proof below. 
\begin{remark}\label{flat}
\rm The structure of linear invariant operators on the locally flat
geometries is well understood in the literature. In particular, it
follows from  the theory of generalized Verma modules that, for each
$k$ as in the theorem, there is exactly one non-trivial operator, up 
to scalar multiples, between the bundles in question, see
e.g. \cite{BoeC}.

It is straightforward to deduce  formulae for these operators: 
First observe that there are preferred scales
in the flat geometries, namely those with $\Rho_{ab}=0$. 
The covariant derivatives commute for such scales
$\xi$, and we will express our formuale in such a scale.
Now for all $0\leq \ell\leq q- 2$
the bundle 
$$
(\yii{\ell}\ce_A) [-\ell] 
$$
is an irreducible component of the $2\ell$-forms on $M$ (appearing with
mutliplicity one) and the operators between the bundles in the theorem
are precisely the non-standard operators in the BGG-resolution of 
the functions,
cf. the diagram in the end of Appendix \ref{A}. Thus it is
clear that the operators concerned are fourth order. 
At the same time, since no primed
indices appear explicitly in our target modules, 
the operators must be given by
$\nd_{ABCD}:=\ep_{A'B'}\ep_{C'D'}\nd_A^{A'}\nd_B^{B'}\nd_C^{C'}\nd_D^{D'}$,
followed by an appropriate $G_0$-module homomorphism onto the target.
(In the preferred scales all curvature vanishes so there is no
possibilty of adding lower order terms.)
Since the covariant derivatives commute
there is only one non-trivial way to apply the Young
projection $\yiiII$ to the image of the operator
$\nd_{ABCD}$. 
Let us write $\tbox_{ABCD}$ for the composition of
such a Young projection with the operator
$\nd_{ABCD}$, followed by the (again unique up to multiple) projection onto 
the desired target.  It is
clear then, that in such a preferred scale for the flat case, the operators
of the theorem are given explicitly by the operator
$\tbox_{ABCD}$. 
\end{remark}

\begin{proof}{Now we are ready to prove theorem \ref{weaks}.} Since 
$\Box_{\alpha\beta\gamma\delta}$ is invariant
we have only  to demonstrate the
claim of the second part of the theorem, namely that for each $k$ as
in the theorem and upon restriction 
to the subbundle 
$$
(\yii{k}\ch_\sigma)[-k]\cong(\yii{k}\ce_S)[-k]
$$
of 
$$
(\yii{k}\ce_\sigma)[-k]
$$
the invariant operator $\Box_{\alpha\beta\gamma\delta}$
takes values in the subbundle 
$$(\wyii{k+2}\ch_\alpha)[-k-2]\cong (\wyii{k+2}\ce_A)[-k-2]$$
of
$(\wyii{k+2}\ce_\alpha)[-k-2]$.  On the way we shall also prove that,
upon restriction to the flat structures, the resulting operator
coincides with the known invariant operator on homogeneous structures
and so it is non-trivial and fourth order. The combination of these
results establishes the theorem.

First we will do the whole task for flat AG-structures. For this case let
us restrict to a scale $\xi$ such that $\Rho_{ab}=0$
on $ M$.  It follows immediately from the definition of
$\Box_{\alpha\beta\gamma\delta}$ in terms of $D^\rho_\alpha$ and the
definition of the latter in terms of $\nd^{R'}_A$ (see \nn{dform}) that
the injecting part of
\begin{equation}\label{optemp}
\Box_{\alpha\beta\gamma\delta}:(\yii{k}\ch_\rho)[w]\to 
(\wyii{k+2}\ce_{\alpha})[w-2]
\end{equation}
is a fourth order operator which, up to a constant non-zero scale, is
a composition of 
$$ \nabla_{ABCD}:  (\yii{k} \ce_R)[w]\to
(\ce_{ABCD})\otimes (\yii{k} \ce_R)[w-2] $$ 
with a Young projection
$$ (\ce_{ABCD})\otimes (\yii{k} \ce_R)[w-2]\to
(\wyii{k+2}\ce_{A})[w-2].  
$$
Since $ k+2\leq q$ the symmetries enjoyed by this are precisely the
symmetries of \nn{optemp} if one formally identifies the twistor
indices of \nn{optemp} with the upper case Roman indices in the
obvious way. Now according to the comments in the remark above (and
given our choice of scale $\xi$), up to scale, all such Young
projections yield the same {\em non-trivial} fourth order operators.
In particular, the injecting part of the image
\begin{equation}\label{op}
\Box_{\alpha\beta\gamma\delta}\left((\yii{k}\ch_\rho)[-k]\right)
\end{equation}
is a non-zero scalar
multiple of the invariant operator $\Box_{ABCD}$ in flat AG-structures. 
That is for 
$$
f_{\rho\cdots \sigma}\in\Gamma(\yii{k}\ch_\rho)[-k]) , 
$$ 
we have that
$$
 \lambda^\alpha_A\cdots \lambda^\gamma_D\lambda^\rho_R\cdots
\lambda^\sigma_S\Box_{\alpha\cdots \gamma} f_{\rho\cdots \sigma} 
$$ 
is independent of the choice of scale $ \xi$ (which recall determines
$\lambda^\beta_A$) from within the preferred class of scales that have
$\Rho_{ab}=0$. It follows that
\begin{multline*}
f_{\rho\cdots \sigma}\mapsto\\
 (\Box_{\alpha\cdots \gamma}f_{\rho\cdots
\sigma}- (Y^A_\alpha \cdots Y^D_\gamma Y^R_\rho\cdots Y^S_\sigma
)(\lambda^{\alpha'}_A\cdots \lambda^{\gamma'}_D\lambda^{\rho'}_R\cdots
\lambda^{\sigma'}_S)\Box_{{\alpha'}\cdots {\gamma'}} f_{{\rho'}\cdots
{\sigma'}} 
\end{multline*}
gives an invariant operator 
$$ 
(\yii{k}\ch_\rho)[-k]\to
(\wyii{k+2}\ce_\alpha)[-k-2] .
$$ 
It is immediately clear that this invariant operator is annihilated
upon contraction with $\lambda^\alpha_A\cdots \lambda^\sigma_S$, where
these $\lambda^\alpha_A$'s are determined by any scale  $\xi'$ from the
preferred class. (That is we do not need $\xi'=\xi$ as the operator is
independent of the choice of scale.) Thus the operator
 vanishes when composed with any such projection
onto the first composition factor of
$$ (\wyii{k+2}\ce_\alpha)[-k-2] .  $$ 
It follows immediately from theorem
\ref{subcomp} 
of appendix \ref{compseries} that the
operator itself must vanish and so the theorem is established for flat 
structures.

It is clear from this result for the flat case, that on general (or 
curved) structures the principal
part of the operator $\Box_{\alpha\beta\gamma\delta}$ on 
$$
(\yii{k}\ch_\rho)[-k]
$$
is non-vanishing and has image in the required bundle.  Now let us fix
a point $q\in M$ and a normal scale $\xi_q$ (see Appendix \ref{normal}) and
consider, at $q$, the composition of this operator with a projection to
an irreducible part of the second composition factor. Notice that our choice
excludes all occurrences of symmetrized 
derivatives of the Rho-tensors (that is the $S$-tensors), since these
vanish under our choices. A typical result
is given by
\begin{equation}\label{trnot}
X^\alpha_{A'}\lambda^\beta_B\cdots \lambda^\sigma_S 
\Box_{\alpha\beta\gamma\delta}f_{\rho\cdots\sigma} . 
\end{equation}
Such a part of the operator must vanish in the flat case and so can
only involve the curvature and its covariant derivatives contracted
into covariant derivatives of the section $f_{R\cdots S}$.  The
unprimed indices of this carry a Young symmetry of the type
$$
\wyiimI .
$$
Now, recall we are considering only torsion-free AG-structures. Thus,
as discussed in Appendix \ref{A}, the only non-zero irreducible component in
the ${\frak g}_0$-part of the curvature $W$ of the normal Cartan connection is
the completely trace-free spinor
  $W^{A'B'D}_{ABC}=\tilde{U}^{A'B'D}_{ABC}=
  \tilde{U}^{[A'B']D}_{(ABC)}$ as the other parts vanish.
This is equivalent to    
$W_{AB\cdots F}:=W^{A'B'D}_{ABC}\vol_{A'B'}\vol_{DE\cdots F}$ which we
will call the
Weyl spinor. Observe that this
has a Young symmetry 
$$
\weylsymm .
$$
The ${\frak g}_1$-part of the curvature $W$ may be expressed
polynomially and purely in terms of the first derivatives of the
latter Weyl spinor (see Appendix \ref{A}). 
Now, from order considerations and classical invariant 
theory it is clear that 
the typical term \nn{trnot} must be a linear 
combination of contractions of the terms 
\begin{equation}\label{terms}
(\nd^{A'}_{A}W_{BC\cdots E})f_{G\cdots H}~~\mbox{ and } 
W_{BC\cdots E}\nd^{A'}_{A}f_{G\cdots H}.
\end{equation}
Considering only the unprimed indices, these terms take values in 
representations of SL$(q)$ described by the tensor product of Young tableaux
$$
\weylsymm \otimes \mb \otimes \yii{k}.
$$
However, we claim that the diagram 
\begin{equation}\label{zero}
\wyiimI{k+2}
\end{equation}
cannot turn up in this tensor product.
To see this note that the only way that one could arrive at the diagram 
\nn{zero} by adding boxes to the diagram 
$$
\weylsymm
$$
is by first producing two full columns, then a further $ 2k+1$ boxes
in appropriate positions. Finally further full columns could be
added. But, since $ q> 2$, for any non-negative integer $ \ell$,
$q+2+2k+1\neq 2q+2k+3+\ell q$ so the outcome is impossible.

Thus the
part \nn{trnot} of the operator must vanish and, by the same argument,
all irreducible parts of second composition factor (i.e. one away from
the injecting part) must vanish.  Thus, by the result (of appendix
\ref{compseries}) that in any composition series \nn{comp}
$V_t=0\implies V_{t+1}=0$ combined with corollary
\ref{strsub}, it follows that the  the operator \nn{op} must take values in
the first composition factor in the bundle
$$
(\wyii{k+2}\ce_\alpha)[-k-2]
$$
 and the theorem is proved.
\end{proof}

\section{Further Observations and Remarks}\label{observations}

As mentioned above local twistors for 4-dimensional conformal spin
structures have been described and investigated by Penrose and others
\cite{Dighton,PenMac,nt}. Analogous local twistor bundles for
complex AG-structures were defined by Bailey and Eastwood in \cite{BaiE}.  The
key to our progress here is the twistor-D operator of definition
\ref{Defn}. This enables a `differentiation' which acts between local
twistor bundles.  Although this operator is new, it is very closely
related to an operator $D_{AP}$ between the so called tractor bundles
of conformal geometry as described in \cite{Goconf} and \cite{Gosrni}.
Much of the calculus surrounding the tractor bundles goes back to Tracy Thomas
whose ideas were recovered and extended in \cite{BaiEGo}. We will not
elaborate in detail on these connections in the current work. However
we briefly indicate here how the twistor-D operator may be used to
define a {\em tractor-D operator} for AG-structures which agrees with
the usual tractor-D operator, as described in \cite{BaiEGo}, on
4-dimensional conformal spin geometries. 

\smallskip\noindent{\bf The Tractor Calculus.}
Let us recall the natural bundles 
$\ce^\alpha\supset\cf^\alpha\simeq \ce^{A'}$. Thus there is the tautological
object $X^{\rho\cdots\sigma}$ providing the identification of the top degree
exterior product of $\cf^\alpha$ with a line bundle: 
\begin{equation}\label{Xdef}
\cf\super{[\rho \cdots \sigma]}{p}=X^{\rho \cdots \sigma}\ce[-1] .
\end{equation}
(In fact $X^{\rho\cdots\sigma}=X^\rho_{R'}\cdots X^\sigma_{S'}\vol^{R'\cdots
S'}$.)
We define a {\em tractor-D} operator, $D_{\alpha\cdots
\beta}$, as follows,
$$
X^{\rho \cdots \sigma}D_{\alpha\cdots \beta} f := 
\Do\super{[\rho \cdots \sigma]}{p}_{\mbox{\tiny $\alpha\cdots \beta$}} f ,
$$
for $f$ (with indices suppressed) in $\ce^{\mu\cdots
  \nu}_{\gamma\cdots \delta}[w]$. Thus, for example, the tractor-D
maps $\ce[w]$ into completely skew valence $p$ cotwistors of weight
$w-1$, $\ce_{[\alpha\cdots \beta]}[w-1]$. Let us call
$\ce_{[\alpha\cdots \beta]}$ the {\em cotractor bundle}. We will use
upper case Greek indices to indicate the abstract indices of the
cotractor bundle and its tensor products and so forth. Thus, for
example, we write
$$
\ce_\Theta=\ce\sub{[\alpha\cdots \beta]}{p} ,
$$
and similarly $\ce^{\Theta}$ for the dual {\em tractor bundle}. 

The tractors and cotractors come from $G$-modules, so they are special cases
of what we have called twistors above. In contrast to the fundamental
twistors, their filtrations are of length $p+1$. We shall see in a moment,
that we recover the tractors of the conformal Riemmanian geometries in the
case $p=2=q$.
 
\smallskip\noindent{\bf The $p=2$ case:}
In this case the tractor bundle is $\wedge^2\ce^\alpha$ and we have 
$$
\ce^{[\alpha\beta]}=\ce^{[AB]}+\ce^{AB'}+\ce^{[A'B']}.
$$
Using the canonical volume form $\vol^{A'B'}$, this may be rewritten
as  
$$
\ce^{\Theta}=\ce^{[AB]}+\ce^{A}_{B'}[-1]+\ce[-1].
$$
In this tractor notation $X^{\rho\sigma}$ of \nn{Xdef} is the canonical
weight one tractor giving the injection $\ce[-1]\to \ce^\Theta$ by
$f\mapsto fX^\Theta $.
In any choice of scale, we have 
$$
X^\Theta=\left( \begin{array}{c} 0\\
                                 0\\
                                 1 \end{array}\right) .
$$
In the cases $q> 2$ no such simplification is available for the analogous
canonical object $Y^{AB}_\Theta$ which describes the injecting part of the
cotractors. Nevertheless it is worthwhile noting that, 
in each choice of scale,
it is given
$Y^{CD}_\Theta=(\delta^C_{[A}\delta^D_{B]},~0~,~0)$. (Here, as above, we
write the injecting part on the left here for consistency with \cite{BaiEGo}.)

Observe that in the $q=2$ case we completely recover the tractors from
\cite{BaiEGo}. In particular,
$h_{\Theta\Lambda}=h_{\alpha\beta\gamma\delta}$ 
is precisely the tractor
metric described in \cite{BaiEGo,Goconf} and in this case
$Y^{AB}_\Theta=X_\Theta \vol^{AB}$ where 
$X_\Theta:=h_{\Theta\Lambda} X^\Lambda$. 

Using the expansions of $\Do^{\rho\sigma}_{\alpha\beta} f$ as in
section \ref{invop}, or
otherwise it is easy to describe explicitly the form of the tractor-D
operator for the $ p=2$ structures. 
Let $\tilde{D}_\Theta$ be the differential operator which, 
in a given choice of scale, may be written $\tilde{D}_\Theta f= (0,\nd_a f
,~wf)$ for $f$ any weight $ w$ twistor (remember that tractor bundles may be
thought of as twistor bundles). This is not itself invariant but in
terms of this the
invariant operator $D_\Theta$ is  given 
$$
D_\Theta f= (w+1)\tilde{D}_\Theta f+Y^{AB}_\Theta \Box_{AB} f
$$
where, again, $f$ is any tractor of weight $w$ and $\Box_{AB}$ is the
operator given by the formula \nn{boxab} (of course $\Box_{AB}$ is only
invariant when $w=-1$). It is easily verified that when $q=2$ this
agrees with the usual formula for the tractor-D operator (apart from
an overall factor of 2, -- compare for
example the formulae in \cite{Gosrni}).  

\smallskip\noindent{\bf Salamon's complex.}  A subcomplex in the
De~Rham complex on a quaternionic manifold $M$ was discussed in
\cite{Sal}. It is just a matter of observation that such a subcomplex
appears for all torsion-free AG-structures with $2=p<q$. This occurs
in the BGG resolution of the sheaf of constant functions, see Figure 1
in Appendix \ref{A} describing the special case $p=2$, $q=4$. Observe
that in that case we can obtain a longer complex if we bypass the bundle in
the vertex of the triangle in Figure 1 
via the second order operator indicated by
the vertical arrow and then continue on the border of the triangle
down to the top degree forms. All this follows immediately from the
fact that the whole diagram, viewed row after row is a genuine
resolution. In fact it is easily verified that this result is typical
and there is an analogous lengthening of Salamon's subcomplex  
for all torsion-free AG-structures with $2=p<q$. 
Using any scale, all the first order operators are always  given
by the appropriate projections of the exterior derivatives expressed
in terms of covariant derivatives.  The `bridging' 
second order operator is given in general by
$$
u^{(A'\cdots C')}_{[A\cdots C]}\mapsto \nd^{S'}_{S}\nd^{(R'}_{R}u^{A'\cdots
C')}_{A\dots C}\ep^{RA\cdots C}\ep_{S'R'}.
$$

\renewcommand{\thesection}{\Alph{section}}
\setcounter{section}{0}
\section{The Cartan connections of AG-structures}\label{A}

\newcommand{\bc}{{\Bbb C}}
\newcommand{\gsl}{{\frak s}{\frak l}}

The AG-structures are specific examples of the so called Cartan geometries.
In general, we have in mind certain 
deformations of homogeneous spaces $G/P$ and
the main defining objects are the Cartan connections on principal
$P$-bundles $\cg$. 
See \cite{Sha} for a complete exposition of the general ideas.

The aim of this Appendix is to apply the general theory to the AG-structures 
and to provide some
background for the main development in this article.

The {\em Cartan connections}
are right invariant forms in $\Omega^1(\cg,{\frak g})$  which reproduce the
fundamental vector fields for the principal action of $P$, and provide
isomorphisms $T_u\cg\to \frak g$ for all $u\in {\cal G}$. The homogeneous
cases are then just the left Maurer-Cartan forms $\omega$ on $G\to G/P$. An
important class among such structures is characterised by two requirements:
the semi-simplicity of $G$, and the existence of the grading of the Lie
algebra ${\frak g}={\frak g}_{-k}\oplus\dots\oplus{\frak g}_k$, $k\in {\Bbb Z}$,
with ${\frak p}={\frak g}_0\oplus\dots\oplus{\frak g}_k$ (the so called
$|k|$-graded Lie algebras). The Lie subgroup $P$ corresponds then to the
subalgebra ${\frak p}$ and it is always a semidirect product of its
reductive part $G_0$ (with Lie algebra ${\frak g}_0$) and the nilpotent
exponential image $P_+$ of ${\frak g}_1\oplus\cdots\oplus{\frak g}_k$. 
In all these cases, the corresponding geometries
are defined in a way similar to classical G-structures and the canonical
bundles ${\cal G}$, together with the canonical Cartan connections, are
constructed from such data. The obstruction against the local equivalence to
the homogeneous spaces is given by the curvature of the Cartan connection,
the two-form $\kappa\in\Omega^2(\cg,{\frak g})$ defined by structure equation 
$$
d\omega =-\frac12[\omega,\omega] + \kappa.
$$ 
By definition, the curvature $\kappa$ is a horizontal two-form and the
presence of the absolute parallelism $\omega$ itself enables us to view
$\kappa$ as the $P$-equivariant function 
$$\kappa:\cg\to {\frak g}_{-}^*\wedge {\frak g}_{-}^*\otimes {\frak
  g}$$
where ${\frak g}_-={\frak g}_{-k}\oplus\cdots\oplus{\frak
  g_{-1}} $ is identified with ${\frak g}/{\frak p}$.  In our case,
the algebra is $|1|$-graded and so the curvature splits into
components $\kappa_{-1}$ (the {\em torsion part}), $\kappa_0$ (the
{\em Weyl part}) and $\kappa_1$.

The canonical Cartan connections are normalised to
have co-closed curvatures $\kappa$, i.e. $\partial ^*\circ \kappa=0$, with
respect to the adjoint to the Lie algebra cohomology differential
$\partial$. 
Such Cartan connections are constructed (including the bundle $\cg$) 
from simple geometric data on the underlying manifold,
see e.g.~\cite{CSch} or \cite{Ta} for explicit constructions in the most
general situations. A very
detailed exposition is also available in \cite{Yam}.

The best known examples are the conformal Riemannian structures and the 
projective geometries, and
all $|1|$-graded cases behave very much similar to them, cf. \cite{Ba,
CSS1, CSS2}. The name
{\em AG-structures} refers in general to
all $|1|$-graded cases where the complexification of  $\frak g$ is ${\frak
s\frak l}(p+q,{\Bbb C})$. In fact, there are only four
relevant series of geometric structures, cf. \cite{KN}:
\begin{itemize}
\item[(1)] ${\frak g}={\frak s\frak l}(p+q,\bc)$ and ${\frak g}_0=
\gsl(p,\bc)\oplus\gsl(q,\bc)\oplus\bc$, ${\frak g}_1=\bc^{q*}\otimes_{\bc}\bc^p$
\item[(2)] ${\frak g}={\frak s l}(p+q,\br)$ and ${\frak g}_0=
\gsl(p,\br)\oplus\gsl(q,\br)\oplus\br$, ${\frak g}_1=\br^{q*}\otimes_{\Bbb R}\br^p$
\item[(3)] ${\frak g}={\frak s l}(p+q,{\Bbb H})$ and ${\frak g}_0=
\gsl(p,{\Bbb H})\oplus\gsl(q,{\Bbb H})\oplus\br$, 
 ${\frak g}_1={\Bbb H}^{q*}\otimes_{{\Bbb H}}{\Bbb H}^p$
\item[(4)] ${\frak g}={\frak su}(p,p)$ and ${\frak g}_0={\frak c
sl}(p,\bc)$, ${\frak g}_1=({\frak su}(p))^*$
\end{itemize}

A general calculus for differential geometry of all $|1|$-graded
geometries was developed in \cite{CSS1}, see also \cite{Sl}.
We are going to review briefly some of the general features of this 
and present
explicit formulae for the AG-structures.

The intuitive explanation of what the geometries look like is as follows:
In each case the tangent space is identified with 
the negative part ${\frak g}_{-1}$ 
of the Lie algebra ${\frak g}$, as a $G_0$-module. 
The most natural choice of the 
Lie group $G_0$
with Lie algebra ${\frak g}_0$ is the adjoint group of the ${\frak
g}_0$-module ${\frak g}_{-1}$.
This choice leads to a sort of minimal
data and in all $|1|$-graded cases this amounts to a
classical G-structure, i.e. a reduction of the general linear frame bundle
to the structure group $G_0$. The Cartan bundle $\cg$
and the Cartan connection $\omega$ are then built out of these data.

In the case of ${\frak g}={\frak sl}(p+q,\bc)$
the structure group described above is a quotient 
$\tilde G_0$ of $G_0=S(\GL(p,\bc)\times \GL(q,\bc))$,
where $G_0\to \tilde G_0$ is a $(p+q)$-fold covering. Thus, it is convenient
to work with the whole $G_0$ instead which, of course, adds some global
structure to our geometries. It does not play any important role locally
though. (In fact, the situation is similar to the spin structures on
conformal Riemannian structures, cf. the case $p=q=2$.) In this paper, 
we are always assuming that this additional structure is given.
Then the $G_0$ structure yields 
an identification of the tangent space of the complex manifold $M$ 
with the tensor product of two auxiliary (complex) vector 
bundles $TM=\ce^A\otimes \ce_{A'}$, together with the fixed isomorphism of
their top degree exterior products, cf. \cite{BaiE}. 

The real split form ${\frak g}={\frak sl}(p+q,\br)$ leads exactly to the 
same description, except we replace complex manifolds and vector bundles by
the real ones, and the reductive group $G_0=S(\GL(p,\br)\times \GL(q,\br))$
equals the minimal structure group $\tilde G_0$ if $p+q$ is odd, while
$G_0\to \tilde G_0$ is a two-fold covering if $p+q$ is even.

The other two real forms are more interesting and quite different, but
we can still include them into the above framework if we deal with the
complex $P$-modules and the complexified tangent bundle
$TM\times_{\Bbb R}\bc$.  Thus we are using the same abstract index
formalism for all these structures, but we have to keep in mind that
it is, with $p$ and $q$ even, the quaternionic form ${\frak
  sl}(\sfrac{p}{2}+\sfrac{q}{2},{\Bbb H})$ which corresponds then to
the discussion of the cases with ${\frak g}={\frak sl}(p+q,{\Bbb
  R})$. This is also compatible with the developments in
\cite{BaiE}, \cite{Sal}. 

Let $\cg$ be the Cartan bundle equipped with the normal connection $\omega$.
The quotient bundle
$\cg_0=\cg/\operatorname{exp}{\frak g}_1$ is a principal fibre
bundle with structure group $G_0$.
Moreover, there is the family of global $G_0$-equivariant sections
$\sigma:\cg_0\to \cg$ parameterised by one forms on $M$ and each such
section $\sigma$ induces the linear connection 
$\gamma^\sigma:=\sigma^*\omega_0$ on $M$ 
(viewed as a principal connection on $\cg_0$). The latter connection,
together with the soldering form $\theta:=\sigma^*\omega_{-1}$ on 
$\cg_0$, forms a Cartan
connection in $\Omega^1(\cg_0,{\frak g}_{-1}\oplus {\frak g}_0)$, and
there is 
the $\sigma$-related Cartan connection $\omega^\sigma\in \Omega^1(\cg,{\frak
g})$. The ${\frak g}_1$-component of the latter connection $\omega^\sigma$ has
to vanish on $T\sigma(T\cg_0)$, while the ${\frak g}_{-1}\oplus{\frak
g}_0$-components of $\omega$ and $\omega^\sigma$ coincide. This implies
that these Cartan connections are related by 
\begin{equation}\label{A-rho}
\omega^\sigma=\omega-\Rho\circ\omega_{-1},
\end{equation}
where $\Rho:\cg\to {\frak g}_{-1}^*\otimes{\frak g}_1$ enjoys the
equivariance properties of a 2-tensor on $M$. The latter tensor is
called the {\em Rho-tensor} defined by the choice of $\sigma$.  The
whole torsion part of the curvature $\kappa$ of the Cartan connection
$\omega$ is constant on the fibres of $\cg$ and provides exactly the
torsion shared by all connections $\gamma^\sigma$.

The absolute parallelism $\omega$ defines the {\em horizontal vector fields}
$\omega^{-1}(X)$ for all $X\in {\frak g}_{-1}$. Now, for each $P$-module $V$
we have the natural vector bundles $\cv$ associated to $\cg$ and their
sections may be viewed as $P$-equivariant functions $s:\cg\to V$. The {\em
invariant differential} $\nd^\omega$ given by the Cartan connection $\omega$ 
is then the obvious differentiation in the directions of the horizontal
vector fields:
$$
\nd^\omega:C^\infty(\cg, V) \to C^\infty(\cg, {\frak g}_{-1}^*\otimes V),
\ \nd^\omega_X s (u) = \omega^{-1}(X)(u).s
$$
In particular, in terms of these invariant derivatives
the Ricci and Bianchi identities have the form
\begin{gather}\label{Ricci}
(\nabla^\omega_X\circ\nabla^\omega_Y - \nabla^\omega_Y\circ\nabla^\omega_X)s =
\lambda(\kappa_{\frak p}(X,Y))\circ s - \nabla^\omega_{\kappa_{-1}(X,Y)}s\\
\label{Bianchi}\sum_{\text{cycl}}\bigl([\kappa(X,Y),Z]-
\kappa(\kappa_{-}(X,Y),Z)-\nabla^\omega_Z\kappa(X,Y)
\bigr)=0
\end{gather}
where $\lambda$ means the representation of $\frak p$ in ${\frak gl}(V)$,
$X,Y,Z\in {\frak g}_{-1}$. 

For irreducible $P$-modules $V$ (and all those with trivial actions of
${\frak g}_1$) we can easily compare the invariant differentials with the
covariant derivatives with respect to any section $\sigma$. We obtain
\begin{equation}\label{1stder}
(\nabla^\omega_X - \nabla^{\gamma^\sigma}_X)s(u)=\lambda([X,\tau(u)])\circ s(u)
\end{equation}
where $\tau:\cg\to {\frak g}_1$ is defined by
$u=\sigma(p(u))\operatorname{exp}\tau(u)$ and it measures the distance of $u$
from the image $\sigma(\cg_0)$ in $\cg$. Consequently, the transformation of
the first derivatives in terms of the change of the scale is
\begin{equation}\label{1sttransf} 
\nabla^{\hat\gamma}_Xs=\nabla^{\gamma}_Xs+\lambda([X,\up])\circ s
\end{equation}

Let us work out this formula in our abstract index formalism. First of all
we need formulae for brackets of elements in ${\frak g}$. We shall write
typical elements $X\in{\frak g}_{-1}$, $Y\in{\frak g}_0$, and $Z\in {\frak
g}_1$ as
$$X=v^A_{A'},\
Y=(u^{A'}_{B'}\delta^B_A + u^B_A\delta^{A'}_{B'}),\
Z=w^{A'}_A
.$$ 
Notice that the convention for ${\frak g}_0$ follows the obvious embedding of
${\frak g}_0$ into the endomorphisms ${\frak g}_{-1}^*\otimes{\frak
g}_{-1}$. In this notation, the brackets in the matrix Lie algebra ${\frak g}$
can be expressed by 
\begin{align*}
[Y, X] &=
- u^{B'}_{A'}v^A_{B'} + u^A_Bv^B_{A'}
\\
[Y, Z] &=
u^{A'}_{B'}w^{B'}_A - u^B_Aw^{A'}_B
\\
[X, Z] &= - w^{A'}_Cv^C_{B'}\delta^B_A +
v^B_{C'}w^{C'}_A\delta^{A'}_{B'}
\end{align*}

Now, the expression ${\frak g}_{-1}\ni X\mapsto [X,\up]\in {\frak g}_0$ 
with  $X=v^A_{A'}$ and $\up=\up^{A'}_B\in{\frak g}_1$, appearing in
(\ref{1stder}), can be be understood as 
$$
v^A_{A'}\mapsto
(-\up^{D'}_A\delta^{A'}_{C'}\delta^C_D+\up^{A'}_C\delta^D_A\delta^{C'}_{D'})
v^A_{A'}
.$$
Thus in order to obtain the formula (\ref{1stder}) we have to act by the 
element
$(-\up^{D'}_A\delta^{A'}_{C'}\delta^C_D+\up^{A'}_C\delta^D_A\delta^{C'}_{D'})$,
viewed as a ${\frak g}_0$-valued one-form with free indices ${}^{A'}_A$, 
composed with the representation $\lambda$. This yields immediately
the formulae in (\ref{ndtrans}).

The Cartan connection $\omega$ induces a connection on all natural bundles 
coming from $G$-modules and the corresponding covariant derivative $\nabla$
is compared to the invariant derivative (and covariant derivatives with
respect to the linear connections $\gamma ^\sigma$) by the formula
\begin{equation}\label{1stcartan}
\begin{align*}
\nabla_X s &= \nabla^\omega_X s +\lambda(X)\circ s 
\\
&= \nabla^{\gamma^\sigma}_X s -\lambda(\Rho.X)\circ s + \lambda(X)\circ s
\end{align*}
\end{equation}
Again, the explicit formulae (\ref{twconn}), (\ref{twconn2}) follow
immediately.

The transformation rule for $\Rho$ under
the change given by $\up$ is then
\begin{equation}\label{rho}
\hat{\Rho}.X=\Rho.X - \nabla_X\up-\tfrac12[\up,[\up,X]]
\end{equation}
In our index formalism this yields exactly (\ref{Ptrans}). 

Next, let us discuss the normalising conditions on the curvatures. The
general formula for the Lie algebra cohomology codifferential $\partial ^*$
(applied to two-forms in ${\frak g}_{-1}^*\wedge {\frak g}_{-1}^*\otimes W$
for a ${\frak g}$-module $W$) reads
$$
\partial^*(Z_1\wedge Z_2\otimes v) = - Z_2\otimes Z_1.v + Z_1\otimes Z_2.v
$$
and so its evaluation on the torsion $T_{ab}{}^c=F^{A'B'C}_{ABC'}+
\tilde F^{A'B'C}_{ABC'}$ where $F^{A'B'C}_{ABC'}=
F^{[A'B']C}_{(AB)C'}$ and $\tilde F^{A'B'C}_{ABC'}=\tilde F^{(A'B')C}_{[AB]C'}$
yields 
$$
\partial^* (T_{ab}{}^c )= 2(-T^{D'B'C}_{ABD'}\delta^{A'}_{C'} +
T^{A'B'D}_{DBC'}\delta^C_A ).
$$
The vanishing of this expression is equivalent to the vanishing of all
traces of the objects $F^{[A'B']C}_{(AB)C'}$, 
$\tilde F^{(A'B')C}_{[AB]C'}$.

Similarly, the evaluation of the codifferential on the ${\frak
g}_0$-component $U^{A'B'C'}_{ABD'}\delta^D_C +
\tilde U^{A'B'D}_{ABC}\delta^{C'}_{D'}$ of the curvature $\kappa$ yields
$$
\partial^*(\kappa_0)= 2(-U^{D'B'A'}_{ABD'} + \tilde U^{A'B'D}_{DBA})
$$
and the condition $\partial ^*\kappa_0=0$ is equivalent to the vanishing of
the two contractions on the right hand side.

By the construction and the general theory, 
the curvatures $\kappa^\sigma$ of the Cartan 
connections $\omega^\sigma$ are $\sigma$-related
to the sum of torsions and curvatures of the induced linear connections
$\gamma^\sigma$ on $\cg_0$. At the same time, the 
relation between $\kappa^\sigma$ and $\kappa$ is
\begin{equation}\label{A-U}
\begin{aligned}
(\kappa^\sigma-\kappa)(u)(X,Y) =\ &\partial \Rho(u)(X,Y) +
\nabla^\omega_X\Rho(u).Y- 
\\
&\nabla^\omega_Y\Rho(u).X + \Rho(u)\circ
\kappa^\sigma_{-1}(u)(X,Y)
.
\end{aligned}
\end{equation}
Our description of the curvature of the twistor connection, see
(\ref{twcurv}), is an immediate consequence of this formula. Furthermore,
the ${\frak g}_0$-component of this expression yields exactly our formula
(\ref{curvdec}). 

The general theory also shows that the whole curvature vanishes if and only
if its harmonic part vanishes and this in turn can be computed explicitly by
the Kostant's version of Bott-Borel-Weil theorem. In our case this means
that the whole curvature is determined by the two components $F$ and $\tilde
F$ of the torsion if $2 < p \le q$. In the case $p=2< q$ only one of the
torsions survives, $\tilde F$, and there appears another invariant component
of $\tilde U^{A'B'D}_{ABC}$, namely the completely trace-free part of $\tilde
U^{[A'B']D}_{(ABC)}$. Let us also notice, that if the torsion happens to
vanish, then the latter component of the Weyl curvature is constant along
the fibres of $\cg\to \cg_0$ and there is no other non-zero component in the
Weyl part of the curvature. 
Moreover, in this case, the ${\frak g}_0$-component of the Bianchi identity
(\ref{Bianchi}) yields for all $X,Y,Z\in {\frak g}_{-1}$
$$
-\partial\kappa_1(X,Y,Z)= \sum_{\operatorname{cycl}}\nabla_Z\kappa_0(X,Y)
.$$
An easy computation reveals that the right hand side is in the kernel of
$\partial$. Because there is no cohomology in that place, the latter
equation has a unique solution for $\kappa_1$ in terms of the derivatives of
the only non-zero component in $\kappa_0$, i.e. of the Weyl spinor $\tilde
U^{[A'B']}_{(ABC)}$. 

\begin{figure}\label{fig1}
\def\uzel#1{\mbox{\small\vbox to0pt{\vss\hbox to0pt{\hss$#1$\hss}\vss}}}
$$
\begin{picture}(280,410)
\put(50,400){\uzel{\ce}}
\put(0,400){\line(1,0){30}}
\put(0,400){\line(0,-1){198}}
\put(0,202){\vector(1,0){15}}
\put(60,390){\vector(1,-1){30}}
\put(100,350){\uzel{\ce^{A'}\otimes\ce_A}}
\put(100,330){\vector(0,-1){60}}
\put(90,340){\vector(-1,-1){30}}
\put(110,340){\vector(1,-1){30}}
\put(50,300){\uzel{(\yiiI\ce_A)[-1]}}
\put(60,290){\vector(1,-1){30}}
\put(-5,300){\line(1,0){20}}
\put(-5,300){\line(0,-1){200}}
\put(-5,100){\vector(1,0){20}}
\put(150,300){\uzel{\yiiI\ce^{A'}\otimes\yiII\ce_A}}
\put(160,290){\vector(1,-1){30}}
\put(140,290){\vector(-1,-1){30}}
\put(150,280){\vector(0,-1){60}}
\put(100,250){\uzel{(\ce^{A'}\otimes\yiIIi\ce_A)[-1]}}
\put(90,240){\vector(-1,-1){30}}
\put(110,240){\vector(1,-1){30}}
\put(100,230){\vector(0,-1){60}}
\put(200,250){\uzel{\yiiiI\ce^{A'}\otimes\yiIII\ce_A}}
\put(190,240){\vector(-1,-1){30}}
\put(210,240){\vector(1,-1){30}}
\put(200,230){\vector(0,-1){60}}
\put(50,200){\uzel{(\yiiII\ce_A)[-2]}}
\put(0,198){\line(1,0){15}}
\put(0,198){\line(0,-1){198}}
\put(0,0){\vector(1,0){25}}
\put(60,190){\vector(1,-1){30}}
\put(150,200){\uzel{\yiiI\ce^{A'}\otimes\yiIIIi\ce_A}}
\put(140,190){\vector(-1,-1){30}}
\put(160,190){\vector(1,-1){30}}
\put(150,180){\vector(0,-1){60}}
\put(250,200){\uzel{(\yiiiiI\ce^{A'})[-1]}}
\put(240,190){\vector(-1,-1){30}}
\put(100,150){\uzel{(\ce^{A'}\otimes\yiIIIiII\ce_A)[-2]}}
\put(90,140){\vector(-1,-1){30}}
\put(110,140){\vector(1,-1){30}}
\put(100,130){\vector(0,-1){60}}
\put(200,150){\uzel{(\yiiiI\ce^{A'}\otimes\ce_A)[-2]}}
\put(190,140){\vector(-1,-1){30}}
\put(50,100){\uzel{(\yiiIII\ce_A)[-3]}}
\put(60,90){\vector(1,-1){30}}
\put(150,100){\uzel{(\yiiI\ce^{A'}\otimes\yiII\ce_A)[-3]}}
\put(140,90){\vector(-1,-1){30}}
\put(100,50){\uzel{(\ce^{A'}\otimes\yiIII\ce_A)[-4]}}
\put(90,40){\vector(-1,-1){30}}
\put(50,0){\uzel{\ce[-6]}}
\end{picture}
$$
\caption{}
\end{figure}

The invariant linear operators between natural bundles over 
locally flat AG-structures are in bijective correspondence with 
the homomorphisms of generalized Verma modules. Thus they are well known from
 representation theory. In particular, all cases with the so called
regular infinitesimal character are obtained by the translation of the
standard De~Rham resolution of the sheaf of constant functions. This is the
source of the celebrated Bernstein-Gelfand-Gelfand resolutions (briefly
BGG resolutions).

The complete BGG resolution of $\ce$ in the special case $p=2$, $q=4$
(i.e. the lowest dimensional interesting quaternionic geometry) is
shown on Figure 1.  The long arrows on the left hand side denote the
non-standard operators.  
One of the aims of our development is to
provide tools for extending such operators to curved geometries. In
fact, there are several methods available, but mostly they fail if
applied to non-standard operators.  
Also the arrows along
the side of the triangle joining $\ce$ and $(\yiiiiI\ce^{A'})[-1]$ 
are worth mentioning. Namely, they 
form the Salamon's subcomplex on quaternionic structures. 

\section{Normal forms for AG-structures}\label{normal}

Given a choice of scale one has a connection $\nd^{\xi}$ 
 on $M$ and for each  point $ q\in M$ one
can define normal coordinates $x^i$ in a neighbourhood of
$q$. Up to a general linear transformation, such coordinates
may be characterized by the conditions that $(1)$, $x^i(q)=0$,  that
 $(2)$ the vectors $\cd/\cd x^i|_q$ give a $G_0$-frame 
at $q$ and that $(3)$ the
coefficients $\Gamma_{(\xi)}$ of $\nd^\xi$, in these coordinates, satisfy
\begin{equation}\label{norm}
\Gamma^i_{jk}x^jx^k=0
\end{equation}
in the neighbourhood where the coordinates are defined. Note that 
$$
\Gamma^i_{(jk)}(q)=0,
$$
and so, at $q$, 
$$
\Gamma^i_{jk}=T_{jk}{}^i .
$$
Similarly differentiating \nn{norm} with respect to the normal
coordinates and evaluating at $q$ we obtain that
$\partial_{(i}\Gamma^i_{jk)}(q)=0$. It follows easily that, at $q\in M$,  
$$
\partial_{k}\Gamma^\ell_{ij}=2R_{k(ij)}{}^\ell +
\frac{1}{6}(3\partial_{k}T_{ij}{}^\ell+2\partial_{(i}T_{j)k}{}^\ell )
+\frac{1}{2}T_{m[k}{}^\ell T_{j]i}{}^m+\frac{1}{2}T_{m[k}{}^\ell T_{i]j}{}^m
.$$
The partial derivatives on the right hand side, of the above, may be
replaced with covariant derivatives at the expense of adding more
terms quadratic in the (undifferentiated) torsion.
By an obvious inductive
argument one can easily continue in this manner and 
recover the following established result.
\begin{proposition}\label{Taylor}
In terms of the normal
  coordinates for $\nd^\xi$, based at $ q\in M$, the coefficients of the
  Taylor series of the $\Gamma_{(\xi)}$ are given by polynomial
  expressions involving the  components of
  the $\nd^\xi$ covariant derivatives of the curvature and torsion
  of 
$\nd^\xi$.
\end{proposition}
Clearly for any choice of scale $\xi$ and $q\in M$ we can
find a $G_0$ family of such normal coordinates.

Fix a choice of scale and normal coordinates in a neighbourhood of
$q\in M$.
Let $u^a$ be a tangent vector at $ q$ and $u^i$ its components in the
normal coordinates. Suppose this is extended to a
section of the tangent bundle in a neighbourhood of $q$ by parallel
transporting $u^a$ along the geodesics through $ q$. Then
$x^i\nd^{\xi}_i u^a=0$ and it is an
elementary   exercise using this to show that the coefficients of 
Taylor series of $ u^a$, about $ q$ and in the normal coordinates, are
given by polynomials in the components  $u^i(q)$ and the
\begin{equation}\label{var}
\Gamma^i_{ij},{}\sub{k\cdots \ell }{t}(q),
\end{equation}
for $t=0,1,\cdots$. (Here $\Gamma^i_{ij},{}_{k\cdots \ell}
:=\cd_\ell\cdots \cd_k \Gamma^i_{ij}$ and the polynomials just
described  are homogeneous of degree 1 in
the components $u^i(q)$.) It follows that the coefficients of the
Taylor series of the {\em normal $G_0$-frame}, corresponding to the
normal coordinates, are polynomial in coefficients $\Gamma^i_{jk}$ and
their normal coordinate derivatives at $ q$. This frame is obtained by
parallel transporting the frame $\cd/\cd x^i|_q$ along the geodesics
through $q$. It follows easily that the coefficients of the connection
$\nd^\xi$ in this normal $ G_0$-frame have normal coordinate Taylor
series with coefficients also polynomial in the variables \nn{var}.
We have a corresponding result for {\em normal spin frames}. These are
constructed as follows. Choose spin frames for $\ce^A(q)$ and
$\ce_{A'}(q)$ consistent with the $G_0$-frame $\cd/\cd x^i|_q$ at $q$
given by the normal coordinates. Now using the spin connections
$\nd^\xi$ parallel transport these frames along the geodesics through
$q$. This determines normal $G_0$-frames for $\ce^A(q)$ and $\ce_{A'}$
in a neighbourhood of $q$. Let $\Gamma^A_{Bi}$ and $\Gamma^{A'}_{B'i}$
be the coefficients of the spin connections with respect to these
frames, where the index $ i$ refers to the normal coordinates (and the
indices $ A,B,A',B'$ here are concrete indices). These coefficients
are linear combinations of the coefficients of the normal $G_0$-frame.
Thus, with the proposition above 
we have the following.
\begin{proposition}\label{normframe} Given a scale $\xi$, and normal
  coordinates at $ x^i$, based at $ q\in M$, let $\Gamma^A_{Bi}$ and
  $\Gamma^{A'}_{B'i}$ be the coefficients of the spin connections with
  respect to the normal spin frame. The coefficients of the Taylor
  series of these functions are given by polynomial expressions
  involving the components of the $\nd^\xi$ covariant derivatives of the
  curvature and torsion of $\nd^\xi$.
\end{proposition}

Given the point $q\in M$ we can can also normalise the scale, at least
formally. Using the equation \nn{Ptrans}, and by considering formal power series, it is easily
verified that one can choose a scale so that 
\begin{equation} \label{normsc}
S\sub{(a \cdots ef)}{s}(q)=0
\end{equation}
for $s=2,3,\cdots, r$ for any given $2\leq r\in {\Bbb N}$.  Let us
suppose that we have chosen and fixed $r$ so that it is sufficiently
large for our calculations and denote this preferred scale $\xi_q$.
\begin{remark}\label{noneedupq} 
In fact it is clear from the form of \nn{Ptrans} that the condition 
\nn{normsc} leaves the 1-jet at $q$ of $\xi_q$ completely free. Thus 
beginning with any scale $\xi$ and an arbitrary point $q\in M$, one
can achieve a normal scale based at $q$, $\xi_q$ by a transformation 
$\xi_q=\Omega \xi$ where $ \up_a(q)=0$. 
\end{remark}

Although we will not use it directly here it is worth observing that, 
in this scale the Taylor series of proposition \ref{Taylor} simplifies
somewhat. Recall the decomposition \nn{curvdec} of the curvature. It
is clear that the jets of the curvature $R^{(\xi)}_{ab}{}^c_d$ are
given linearly by the jets of the
tensor $U^{(\xi)}_{ab}{}^c_d$ and the jets of the Rho-tensor 
$\Rho^{(\xi)}_{ab}$. Considering various Young projectors acting on 
$\nd_a\nd_b\cdots \nd_d \Rho^{(\xi_p)}_{ef}$ one easily concludes 
that, at $q\in M$, this tensor is determined by   
$\nd_a\nd_b\cdots \nd_d \Rho^{(\xi_p)}_{[ef]}$,
$\nd_a\nd_b\cdots \nd_{[d} \Rho^{(\xi_p)}_{e]f}$ and lower order
terms. But, by \nn{skewrho} $\Rho^{(\xi)}_{[ef]}$ is given by a linear
formula in terms
of a $\nd^{\xi}$ derivative of the torsion. Thus we obtain the
following simplification to the above proposition.
\begin{proposition}\label{modTaylor}
  Let $q\in M$ and $\xi_p$ be a scale such that \nn{normsc} is satisfied.
  Let $x^i$ be normal coordinates for $\nd^{\xi_p}$ based at $q$.
  Then, in terms of these coordinates, the coefficients of the Taylor
  series of the $\Gamma_{(\xi_p)}$ (to order $r+1$) are given
  by polynomial expressions in the components of the covariant derivatives of
  the torsion $T^{(\xi_p)}_{ab}{}^c$ of
  $\nd^{\xi_p}$ and the components of the covariant derivatives
 of the tensors $U^{(\xi_p)}_{ab}{}^c_d$ and $\nd_{[a}
  \Rho^{(\xi_p)}_{b]c}$.
\end{proposition}

\section{Composition series}\label{compseries}

\newcommand{\bY}{{\bf Y}}

In the following discussion we will review several notions and terms
for representations of a group $H$. We have, for the most part, not
said anything about the nature of this group since an explicit
description of the group is not required for most of the results here.
Of course for application of these results to the other parts of this
article one may take $H$ to be a parabolic $P$ in one of real Lie
algebras $G$ as discussed in the introduction. We would also like to
point out that in this case the terms introduced (such as
``composition series'' and ``injecting part'' etcetera) can be adapted
in an obvious way to the natural bundles that $P$ induces and indeed to
differential operators that take values in such natural bundles. Throughout
the article we have used this observation without other
mention.

Suppose $V$ is an $H$-module for some group $H$. Let $W$ be an 
$H$-submodule of $V$ then we have an exact sequence 
$$
0\to W\to V\to U\to 0
$$  
where $U$ is the required quotient. Following Buchdahl (see also
\cite{BaiEGo}) it is often convenient to express this
as a composition series in the following schematic manner,
$$
V=U+W .
$$

Suppose now that $V$ is any non-trivial finite dimensional module for
the group $H$.  We construct a composition series of $V$ as follows.
Let $V_1^1$ be an irreducible submodule of $V$. If there is a
non-trivial submodule of $V$ in a complement to $ V_s^1$ then there is
at least one irreducible one which may denote $ V_1^2$. Continuing in
this manner suppose that $\{ V_1^1,V_1^2,\cdots ,V_1^{m_1}\}$ is a
{\em maximal} set of such submodules, meaning that there are no
non-trivial submodules of $V$ in a complement to $V_1:=
\oplus_{i=1}^{m_1} V_1^i$.  We call $ V_1$ the first {\em composition
factor} of $V$, while the irreducible submodules $V_1^i$
($i\in\{1,\cdots, m_1\}$) in this, will be described as {\em injecting
parts} of $V$.

Now let $U_{2}:= V/V_1$. Then $U_{2}$ is an $H$-module and so we
may similarly choose a set of irreducible submodules of this, $
V_{2}^i$, $i=1,\cdots ,m_{2}$, such that this set is maximal in $U_{2}$. 
We write $V_{2}$ for the first composition factor of $ U_{2}$, that is 
$V_{2}=\oplus_{i=1}^{m_{2}}V_{2}^i$. 

Now we may consider $U_{3}:=U_{2}/V_{2}$ and seek a maximal set of
irreducible submodules of this (which we denote $V_{3}^i,
~i=1, \cdots ,m_{3}$) and so on. Note that at any stage $V_{t}=0$
if and only if $U_{t}=0$. Since the $ U_{t+j}$, for $j\geq 1$, are
quotients of $U_{t}$ it follows that $ V_{t}=0$ implies $
V_{t+j}=0$ for all $ j\geq 1$. In fact since $V$ is assumed finite
dimensional it clear that there exists some positive integer $r$ such
that $ V_{r+1}=0$ while $V_{r}\neq 0 $. With that determined the 
{\em composition series} of $V$ is given,
\begin{equation}\label{comp}
V=(\oplus_{i=1}^{m_r}V^i_r) +(\oplus_{i=1}^{m_{r-1}}V^i_{r-1})+\cdots +(\oplus_{i=1}^{m_1}V^i_1) .
\end{equation}
We describe $V_k=\oplus_{i=1}^{m_k}V^i_k$ as the $k^{th}$ {\em
composition factor} of $ V$.  The  $V^i_r$ ($i\in\{1,\cdots,m_r\}$) will 
be called the {\em projecting parts} of $V$.
(It is usual to describe $V_r+V_{r-1}+\cdots + V_1$ as the composition series 
for $ V$. For our purposes it is convenient to choose a decomposition of the 
composition factors as indicated.)

We have the following results.
\begin{theorem}\label{subcomp}
Suppose an $H$-module $V$ has a composition series as in \nn{comp}. Then 
for $S$  an $H$-submodule of $ V$ we have 
$$S\cap V_1=0 \Leftrightarrow S=0 .$$
\end{theorem}
\begin{proof}
  The implication $\Leftarrow$ is clear. Suppose now $S$ is an
  \ul{irreducible} $ H$-submodule such that $S\cap V_1=0$. Then $S=0$ since
  $\{V_1^1,\cdots , V_1^{m_1}\}$ is a maximal set of irreducible
  submodules of $ V$. Now suppose $ S$ is any $ H$-submodule such that
  $S\cap V_1=0$. Then an irreducible $ H$-submodule $S'$ of $S$ is an
  irreducible $ H$-submodule of $ V$ such that $S'\cap V_1=0$. Thus by
  the established result $ S'=0$. Thus $S$ has no non-trivial
  irreducible submodules and so $S=0$ as claimed. 
\end{proof}
The
following indicates that a composition series is unique up to
some possible choice for the splitting of each part into irreducibles.
\begin{corollary}\label{unique}
Suppose an $H$-module $V$ has a composition series as in \nn{comp} and also a 
composition series 
$$V=(\oplus_{i=1}^{\tilde{m}_{\tilde{r}}}\tilde{V}^i_{\tilde{r}} +
(\oplus_{i=1}^{\tilde{m}_{{\tilde{r}}-1}}\tilde{V}^i_{{\tilde{r}}-1})+\cdots 
+(\oplus_{i=1}^{\tilde{m}_{1}}\tilde{V}^i_{1})
$$
then $\tilde{r}=r$, $\tilde{m}_1=m_1,\cdots , \tilde{m}_r=m_r$ and
$V_1=\tilde{V}_1:=\oplus_{i=1}^{m_1}\tilde{V}^i_1, \cdots ,
V_r=\tilde{V}_r:=\oplus_{i=1}^{m_r}\tilde{V}^i_r$. 
Furthermore in each composition factor
$V_k$ one can arrange the numbering of the $V_k^i$ so that for each
$i\in\{1,\cdots ,m_k\}$ $V_k^i\cong \tilde{V}_k^i$. If for any $i$ the module 
$V_k^i$ occurs with multiplicity one in $V_k$ then we get 
$V_k^i= \tilde{V}_k^i$.   
\end{corollary}
\begin{proof}
  The first part of this is immediate by repeated application of the
  theorem while the last part follows from Schur's lemma.
\end{proof}
{}From this in turn we get the following corollary.
\begin{corollary}\label{strsub}
 Let  $V$ be an $H$-module with composition series as in  \nn{comp}. 
  If $S$ is an $H$-submodule of $V$ then $ S$ has a composition series
$$
S=(\oplus_{i=1}^{\ell_{r_s}}S^i_{r_s}) +(\oplus_{i=1}^{\ell_{({r_s}-1)}}S^i_1)
+\cdots +(\oplus_{i=1}^{\ell_1}S^i_{1})
$$ 
where for each $k\in \{1, \cdots , r_s\}$ and $i\in \{1,\cdots ,\ell_k\}$
there is some $j\in \{1,\cdots , m_{k}\}$ such that
$$
S^i_k\cong V^j_{k}, 
$$
with equality if $V^j_{k}$ occurs with multiplicity one in $
V_{k}$.
\end{corollary}
Thus all homomorphisms between finite dimensional $H$-modules
$ V$ and $ W$ are determined by the composition series for $ V$ and $
W$, at least up to an isomorphism ambiguity due to the multiplicity of
irreducible components in each part.

We are in particular interested in the composition series of $ P$
modules which are the restriction to $ P$ of irreducible $ G$ modules
and also their $ P$-submodules. Recall that $P$ is a maximal parabolic
in a group $G$ which is a   real form of the complex 
semisimple groups $\SL(p+q,{\Bbb C})$. 
In this case some aspects of the composition series are rather easily 
described.

Let $ V_\alpha$ be the dual to the standard representation of $ G$.
Then we have an exact sequence of $ P$-modules
\begin{equation}\label{seq}
0\to V_A\to V_\alpha \to V_{A'}\to 0.
\end{equation}
Let $ Y^A_{\alpha}$ be the canonical element of $ V^A\otimes V_\alpha$
giving the injection $V_A\to V_\alpha$ and $ X_{A'}^\alpha$ be the
canonical element of $ V_{A'}\otimes V^\alpha$ giving the surjection $
V_\alpha \to V_{A'}$. (This notation is borrowed from the notation for
the corresponding objects for bundles these modules induce.) Let us also write 
$H_\alpha $ for the image of $ V_A$ in $ V_\alpha$.

Now irreducible $ G$-modules may be described by Young diagrams. Using 
notation as in section \ref{newops} 
we may write for example
$$
\bY(b) V_\alpha ,
$$
where $ \bY(b)$ indicates a Young diagram with a total of $ b$
boxes. (We will suppose the height of this diagram is no greater than
$p+q$ so this module is not trivial.)  Elements of this module consist
of vectors which carry $b$ indices,
$$
v\sub{\alpha\cdots \gamma}{b},
$$
and a symmetry indicated by the Young diagram. Regard this now as a $
P$-module by restriction and consider the subspace of  vectors
that have the property that they are $ X$-{\em saturated}, that is
they are annihilated upon contraction with $X^\alpha_{A'} $ on any
index,
$$
0=X_{A'}^\alpha v_{\alpha\beta \cdots \gamma}
=X_{A'}^\beta v_{\alpha\beta \cdots \gamma}=\cdots
=X_{A'}^\gamma v_{\alpha\beta \cdots \gamma}.
$$
The space of such vector clearly forms a $P$-submodule of $\bY(b) V_\alpha $. 
Considering each index in turn it is clear that it is a submodule of 
$ \otimes^b H_\alpha$. Thus it is precisely the submodule 
$$
\bY(b) H_\alpha .
$$
Of course this may be trivial but in any case
$$
\bY(b) H_\alpha \cong \bY(b) V_A
$$
and so it is irreducible. Thus if this is not zero then it gives
the unique injecting part of $ \bY(b) V_\alpha$ (which is therefore
also the first composition factor). If the height of the diagram $
\bY(b)$ is no greater than $q$ then we are in this situation, that is 
$\bY(b) V_A\neq 0 $, and we shall
henceforth assume this is the case since it is sufficient for our
purposes.  The quotient
$$
(\bY(b) V_\alpha )/V_1
$$
may clearly be identified with the direct sum of the distinct images of
$\bY(b) V_\alpha $ under the mapping given by contraction with one $
X^\rho_{R'}$. 
Each of these distinct images carries a Young symmetry
on its twistor indices (that is the greek indices) and, reasoning
essentially as for the previous case, one sees that the irreducible
parts of the second composition factor are the subspaces of images
that are annihilated by any (further) contraction with $
X^\sigma_{S'}$. One can clearly continue in this manner to determine
the entire composition series. For the purposes of this article we
only explicitly require an understanding of this to the level
described. Let us just finally observe that given a choice of
splitting of the sequence \nn{seq}, 
or equivalently a choice of
$\lambda^\beta_{B}$ such that 
$\lambda^\beta_{B} Y_\beta^A=\delta^A_B$ it follows immediately from
the observations here that we may describe these parts of the
composition series as follows. The injecting part of $ \bY(b) V_\alpha$
may be identified with the space of vectors $
\lambda^\alpha_A\lambda^\beta_B\cdots \lambda^\gamma_C v_{\alpha\beta
\cdots \gamma}$ for $v_{\alpha\beta \cdots \gamma}\in \bY(b)V_{\alpha}
$. The second composition factor may similarly be identified 
with the vector space of objects consisting of vectors in $ \bY(b) V_\alpha$
contracted into $(b-1)$ $ \lambda^\rho_{R}$'s and one $ X^\rho_{R'}$.
The corresponding result for induced bundles is used in section \ref{newops}.

\end{document}